\documentclass[hidelinks,12pt,a4paper]{article} 

\usepackage[english]{babel}
\usepackage{latexsym}
\usepackage{amssymb}
\usepackage{amsmath}
\usepackage{amsfonts}
\usepackage{mathtools}
\usepackage{graphics}
\usepackage{enumitem}
\usepackage{afterpage}
\usepackage{fancybox}
\usepackage{hyperref}
\usepackage{xspace}
\usepackage{xcolor}
\usepackage{color}
\usepackage{scalerel}
\usepackage{tikz}
\usetikzlibrary{graphs}
\usetikzlibrary{arrows}
\usetikzlibrary{matrix}

\newcommand{\al}{\alpha}
\newcommand{\be}{\beta}

\newcommand{\ep}{\varepsilon}
\newcommand{\la}{\lambda}

\newcommand{\ita}{\textit}

\newcommand{\mN}{\mathbb{N}}

\newcommand{\mT}{\mathbb{T}}

\newcommand{\sN}{\mathcal{N}}
\newcommand{\sP}{\mathcal{P}}

\newcommand{\sW}{\mathcal{W}}

\newcommand{\fa}{\mathfrak{a}}
\newcommand{\fb}{\mathfrak{b}}
\newcommand{\fc}{\mathfrak{c}}

\newcommand{\wqo}{\textsf{wqo}\xspace}
\newcommand{\qo}{\textnormal{qo}\xspace}
\newcommand{\po}{\textnormal{po}\xspace}

\newcommand{\imp}{\Rightarrow}
\newcommand{\frec}{\rightarrow}

\newcommand{\x}{\cdot}
\newcommand{\sse}{\leftrightarrow}
\newcommand{\et}{\wedge}
\newcommand{\vel}{\vee}

\newcommand{\tri}{\trianglelefteq }
\newcommand{\tild}{\thicksim}
\newcommand{\subeq}{\sqsubseteq}

\newcommand{\qed}{\hfill$\square$}

\newcommand{\KTlw}{K$\mbox{T}_\ell(\omega)$}
\newcommand{\KTln}{K$\mbox{T}_\ell(n)$}
\newcommand{\KTw}{K$\mbox{T}(\omega)$}
\newcommand{\KTn}{K$\mbox{T}(n)$}

\newcommand{\RCA}{\textsf{RCA$_0$}}

\newcommand{\ACA}{\textsf{ACA$_0$}}

\renewcommand{\preceq}{\preccurlyeq}
\renewcommand{\leq}{\leqslant}
\renewcommand{\geq}{\geqslant}
\renewcommand{\nleq}{\nleqslant}

\renewcommand{\emph}[1]{\textsc{#1}}
\newcommand{\defAs}{ \vcentcolon \Leftrightarrow}
\newcommand{\putAs}{\! \vcentcolon=}

\newcommand{\proof}{\noindent \ita{Proof}\ }

\newtheorem{theor}{Theorem}[section]
\newtheorem{prop}{Proposition}[section]
\newtheorem{lem}{Lemma}[section]
\newtheorem{cor}{Corollary}[section]

\newtheorem{defi}{Definition}[section]
\newtheorem{conj}{Conjecture}[section]

\newtheorem{rem}{Remark}[section]

\newcommand{\Rem}[1]{Remark~\ref{#1}}
\newcommand{\Theor}[1]{Theorem~\ref{#1}}
\newcommand{\Lem}[1]{Lemma~\ref{#1}}
\newcommand{\Def}[1]{Def.~\ref{#1}}

\newcommand{\Sec}[1]{Sec.~\ref{#1}}

\newcommand{\COMMENT}[1]{}

\definecolor{navy}{rgb}{0,0,0.5}
\definecolor{magenta}{rgb}{0.7,0,0.7}

\def\sqdot{\mathbin{\scalerel*{\strut\rule{1ex}{1ex}}{\cdot}}}

\graphicspath{{../grafici/}}

\title{Proof-Theoretic Relations between \\
Higman's and Kruskal's theorem,\\
and Independence Results\\
for Tree-like Structures}     
\author{Gabriele Buriola$^1$, Andreas Weiermann$^2$}
\date{%
    $^1$University of Verona\\ \vspace{0.1 cm}%
    $^2$Ghent University
}

\begin{document}

\maketitle

\begin{abstract}

Higman's lemma and Kruskal's theorem are two of the most celebrated results in the theory of well quasi-orders. In his seminal paper G. Higman obtained what is known as Higman's lemma as a corollary of a more general theorem, dubbed here as Higman's theorem. While the lemma deals with finite sequences over a well quasi-order, the theorem is about abstract operations of arbitrarily high arity. 

J.B. Kruskal was fully aware of this broader framework: in his seminal paper, he not only applied Higman’s lemma at crucial points of his proof but also followed Higman’s proof schema. At the conclusion of the paper, Kruskal noted that Higman’s theorem is a special case—restricted to trees of finite degree—of his own tree theorem. Although he provided no formal reduction, he included a glossary translating concepts between the tree and algebraic settings. The equivalence between these versions was later clarified by D. Schmidt and M. Pouzet. In this work, we revisit that equivalence to illuminate the proof-theoretical relationships between the two theorems within the base system 
\RCA\, of reverse mathematics.

Moreover, some independence results over first- and second-order theories are treated. In particular, tree-like structures, involving either Ackermannian terms or exponential expressions, are studied unveiling well-foundedness properties that are independent from Peano arithmetic and relevant fragments of second-order arithmetic.


\end{abstract}

\bigskip


\smallskip

\noindent{\bf MS Classification 2020:}\ \textbf{03B30}, 03F35, 05C05, 06A07.

\section*{Introduction}\label{sec:introduction}

Higman's lemma and Kruskal's theorem are two of the most celebrated results in the theory of \ita{well quasi-orders}, routinely shortened in {\wqo}. In his seminal paper \cite{Higman52} Higman obtained what is known as Higman's lemma as a corollary of a more general theorem, dubbed here Higman's theorem. While the lemma, which says that if $Q$ is a {\wqo}, then the set $Q^*$ of finite sequences on $Q$ is {\wqo}, deals with finite sequences over a well quasi-order  (and so, implicitly, with only the binary operation of juxtaposition), the theorem is about abstract operations of arbitrary high arity, covering a  far more extensive spectrum of situations. 

Kruskal was well aware of this more general setup; in his time-honoured paper \cite{Kruskal60} not only did he use Higman's lemma in crucial points of the proof of his own theorem, but also followed the same proof schema as Higman. 
Moreover, in the very end of his article, Kruskal explicitly stated how Higman's theorem is a special version of Kruskal's  tree theorem; namely a restriction to trees of finite branching degree, i.e., trees with an upper bound regarding the number of immediate successors of each node. Although no proof of the reduction is provided, he presented a glossary to properly translate concepts from the tree context of his paper to the algebraic context of Higman's work.  The equivalence between restricted versions has subsequently been exposed by Schmidt \cite{Schmidt79},  whereas Pouzet \cite{Pouzet85} has given, together with that equivalence, an infinite version of Higman's theorem  which proves equivalent to the general Kruskal theorem.

Besides its crucial role in theoretical computer science, e.g., for term rewriting \cite{Dershowitz82,Dershowitz87}, 
Kruskal's theorem characterizes some relevant parts of proof theory. In the vein of the Friedman--Simpson programme 
of reverse mathematics \cite{Simpson09}, considerable efforts have been undertaken in order to properly 
understand the proof-theoretic strength of the tree theorem, as well as of some of its most important versions \cite{Freund20gap,Gallier91,RW93}; nevertheless, quite a few aspects remain obscure.

We revisit  the aforementioned equivalences to obtain a clear view of the proof-theoretical relations between the different versions of Higman's and Kruskal's theorem, paying particular attention to the former's algebraic formulation. More precisely, we establish, over the base theory \RCA, the equivalence of the finite and infinite version of Higman's and Kruskal's theorem obtaining, as an immediate corollary, the fact that the finite case of Kruskal's theorem for labelled trees is strictly stronger than the infinite case of the unlabelled version, using our future notation KT$(\omega) \not   \imp \forall n\,$KT$_\ell(n)$. The ultimate goal is to complete the picture within the frame of reverse mathematics and ordinal analysis, following \cite{RW93}. 

On top of that, and extending the findings of \cite{FW23}, some independence results over first- and second-order theories are treated. In particular, tree-like structures, involving either Ackermannian terms or exponential expressions, are studied unveiling well-foundedness properties that are independent from Peano arithmetic and relevant fragments of second-order arithmetic.

\smallskip

The article has the following structure: \Sec{sec:wqo} briefly collocates the theory of {\wqo} in reverse mathematics, mainly over \RCA, and contains some novel analyses regarding closure properties for ${\wqo}$ under order-preserving and order-reflecting functions; in \Sec{sec:algebras} we introduce the two main topics of the paper, Higman's and Kruskal's theorem, as well as their ingredients, namely \ita{ordered algebras} and \ita{trees}; \Sec{sec:HK}, which is the core of the article, is devoted to the proof-theoretical relations, over \RCA, between different versions of Higman's and Kruskal's theorem, proving in particular the equivalence between, respectively, the two finite versions and the two infinite versions; in \Sec{sec:HD} other proof-theoretical results, related this time to Higman's lemma and Dickson's lemma, are exposed; \Sec{sec:ind} treats some independence results over first- and second-order theories (prominently Peano arithmetic and fragments of second-order arithmetic) of properties related to tree-like structures; \Sec{sec:work} consists in a brief presentation of related works and further possible developments; a conclusion section ends the paper.

Part of the material in this article, in particular \Sec{sec:HK} and \Sec{sec:HD}, has been treated in the PhD thesis of the first author \cite{BuriolaPhD}.

\section{Well Quasi-Orders in Reverse Mathematics}\label{sec:wqo}

Given the crucial role played by well quasi-orders in our results, this section is dedicated to their formal introduction in reverse mathematics, particularly over the base theory \RCA, and to a first analysis of their properties which include novel analysis regarding closure results for ${\wqo}$ under order-preserving and order-reflecting functions. For an introduction to reverse mathematics, including the axioms of \RCA, we refer to \cite{Simpson09}; while for a recent introductory survey on {\wqo}, including the more refined notion of better quasi-orders, see~\cite{BS:wqo}.

\smallskip

All the following definitions are given in \RCA.

\begin{defi}\label{def:qoRCA0} A \emph{quasi-order}, {\qo}, $(Q, \leq)$ is given by a subset $|Q| \subseteq \mN$ together with a binary relation $\leq \, \subseteq |Q| \! \times \! |Q|$ which is reflexive and transitive; if in addition $\leq$ is antisymmetric, then $(Q,\leq)$ is a \emph{partial order}, {\po}.
\end{defi}
In the following, $|Q|$ is simply denoted by $Q$, and we may denote with $Q$ also the qo $(Q,\leq)$, omitting the quasi-order relation $\leq$; moreover, we use $ \tild $ for the comparability relation, namely $p \! \tild \! q \defAs p \! \leq \! q  \vel q \! \leq \! p$. Observe that, provable in \RCA, if $(Q,\leq)$ is a quasi-order, than $(Q/\! \! \tild, \leq)$ is a partial order.

To properly work with {\qo}, some auxiliary definitions are needed.
\begin{defi}\label{def:auxRev} For every quasi-order $(Q,\leq)$, 

\begin{itemize}

\item the \emph{closure} of a subset  $ B$ of $Q$ is given by 
$\bar{B} = \{ q \! \in \! Q \, |\, \exists b \! \in \! B \, b \! \leq \! q \}$;

\item a subset of $Q$ is \emph{closed} if it equals its own closure, and a closed subset is \emph{finitely generated} if it is the closure of a finite set;

\item a \emph{sequence} $(q_k)_k$ in $Q$ is a total function from $\mN$ to $Q$;  

\item an \emph{antichain} in $Q$ is a sequence $(q_k)_k$ in $Q$ such that $q_i$ and $q_j$ are incomparable, i.e., $q_i \! \nleq \! q_j \land q_j \! \nleq \! q_i$, whenever $i \neq j$;

\item an \emph{extension} of $(Q, \leq)$ is a {\qo} $\preceq$ on $Q$ extending $\leq$ in the sense that  $p\! \leq \! q\Rightarrow p \! \preceq \! q$ and such that for all $p$ and $q$, $p \! \preceq \! q \et q \! \preceq \! p \Rightarrow p \! \tild \! q$.\footnote{The second condition ensures a bijection between the extensions of the {\qo} $(Q,\leq)$ and the extension of the {\po} $(Q/ \! \! \tild,\leq)$.}

\end{itemize}

\end{defi}

We move on to well quasi-orders.
\begin{defi}\label{def:wqo} A {\qo} $(Q,\leq)$ is a \emph{well quasi-order}, {\wqo}, if for every sequence $(q_k)_k$ in $Q$ there exist indexes $i < j$ such that $ q_i \leq q_j$.
\end{defi}
Many different definitions have been proposed for well quasi-orders; in his classical historical survey \cite{Kruskal72}, J. Kruskal referred to {\wqo} as a ``frequently discovered concept''. It is therefore natural to ask if all these definitions are equivalent in Reverse Mathematics, in particular over weak theories like \RCA \, or \textsf{WKL$_0$}. This problem has been thoroughly treated by Cholak et al. \cite{CMS04} and by Marcone \cite{Marcone05,Marcone20}. For a recent introductory survey on {\wqo} and the stronger notion of better quasi-orders see \cite{BS:wqo}. We now briefly summarize their findings.

Let us consider the following definitions for {\wqo}.
\begin{defi}\label{def:wqosRev} Let $(Q,\leq)$ be a qo, then $Q$ is:

\begin{enumerate}

\item a \emph{sequentially well quasi-order}, {\wqo}(set), if every sequence $(q_k)_k$ in $Q$ has an infinite ascending subsequence, i.e., there are indices $k_0<k_1< \ldots$ such that $q_{k_i} \! \leq \! q_{k_j}$ whenever $i \! < \! j$;

\item an \emph{antichain well quasi-order}, {\wqo}(anti), if $Q$ has no infinite descending chains and no infinite antichains;

\item an \emph{extensional well quasi-order}, {\wqo}(ext), if every linear extension $\preceq $ of $\leq$ is well-founded;

\item {\wqo}(fbp) if $Q$ has the \emph{finite basis property}, i.e., every closed subset is finitely generated.

\end{enumerate}

\end{defi}
For what concerns the relations over \RCA, \textsf{WKL$_0$} and \ACA, between these definitions, the main results in \cite{CMS04,Marcone05,Marcone20} are the following.

\begin{theor}[Cholak, Marcone, Solomon]\label{theor:wqorelation} Consider the previous definitions: {\wqo}, {\wqo}(set), {\wqo}(anti), {\wqo}(ext), {\wqo}(fbp). Then over

\begin{description}

\item[RCA$_0$] {\wqo} is equivalent to {\wqo}(fbp) and all the other relations are exhaustively given by the transitive closure of the following schema, namely no arrow can be inverted:

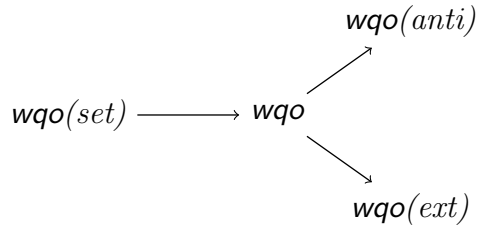
\begin{center}
\begin{tikzpicture}[every node/.style={rectangle, fill=white, text=black}]
\node (wqoset) at (-2.75,0) {{\wqo}(set)};
\node (wqo)    at (0,0)     {{\wqo}};
\node (wqoanti)   at (1.75,1.25)  {{\wqo}(anti)};
\node (wqoext)    at (1.75,-1.25)  {{\wqo}(ext)};

\draw[line width=0.5pt,  ->] (wqoset)--(wqo);
\draw[line width=0.5pt,  ->] (wqo)--(wqoanti);
\draw[line width=0.5pt,  ->] (wqo)--(wqoext);

\end{tikzpicture}

\end{center}

\item[WKL$_0$] The exact implications over WKL$_0$ are given by the transitive closure of the following schema:

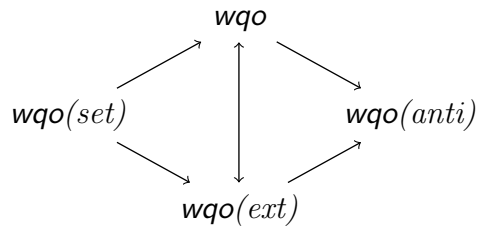
\begin{center}
\begin{tikzpicture}[every node/.style={rectangle, fill=white, text=black}]
\node (wqoset) at (-2.25,0) {{\wqo}(set)};
\node (wqo)    at (0,1.25)     {{\wqo}};
\node (wqoanti)   at (2.25,0)  {{\wqo}(anti)};
\node (wqoext)    at (0,-1.25)  {{\wqo}(ext)};

\draw[line width=0.5pt,  ->] (wqoset)--(wqo);
\draw[line width=0.5pt,  ->] (wqoset)--(wqoext);
\draw[line width=0.5pt,  <->] (wqo)--(wqoext);
\draw[line width=0.5pt,  ->] (wqo)--(wqoanti);
\draw[line width=0.5pt,  ->] (wqoext)--(wqoanti);

\end{tikzpicture}

\end{center}

\item[ACA$_0$] All the definitions are equivalent.

\end{description}

\end{theor}

\proof See \cite{CMS04,Marcone05,Marcone20}. \qed

Given their good behavior, one of the main problems in the theory of well quasi-orders concerns how to obtain new {\wqo}'s or, equivalently, how to preserve the property of being {\wqo}. Two standard operations on {\qo} are \ita{product} and \ita{(disjoint) union}.
\begin{defi}\label{def:wqoOp} Let $(P,\leq_P)$ and $(Q,\leq_Q)$ be {\qo}, then we define:

\begin{enumerate}

\item the \emph{(disjoint) union} $(P \, \dot{\cup} \, Q, \leq_{P \, \dot{\cup} \, Q})$ with
\[
p \leq_{P \, \dot{\cup} \, Q} q \defAs  (p, q \! \in \! P \et p \! \leq_P \! q ) \, \vel \, (p, q \! \in \! Q \et p \! \leq_Q \! q) ;
\]

\item the \emph{product} $(P \times Q, \leq_{P  \times Q})$ with 
\[ (p_1, q_1) \leq_{P \times Q} (p_2, q_2) \defAs p_1 \! \leq_P \! p_2 \, \et \, q_1 \! \leq_Q \! q_2.
\]

\end{enumerate}

\end{defi}
The word ``disjoint'' may be omitted; for other possible quasi-order operations, e.g., the \ita{sum} $P +Q$, see \cite{BW:ordinals} and \cite{Marcone20}.

We now consider the closure of {\wqo} under product, union and subsets, starting with the good behavior of union.

\begin{lem}[Marcone]\label{lem:wqounion} Let $\sP$ be any of the properties {\wqo}, {\wqo}(set), {\wqo}(anti) or {\wqo}(ext); if $(P,\leq_P)$ and $(Q,\leq_Q)$ satisfy property $\sP$, then \RCA\, suffices to prove that $(P \, \dot{\cup} \, Q, \leq_{P \, \dot{\cup} \, Q} )$ has property $\sP$.
\end{lem}

\proof See \cite{Marcone20}. \qed

Except for {\wqo}(ext), a similar well behavior also holds for subsets.
\begin{lem}\label{lem:wqosub} Let $\sP$ be any of the properties {\wqo}, {\wqo}(set), {\wqo}(anti); if {$(Q,\leq_Q)$} satisfies property $\sP$ and $P \! \subseteq \! Q$, then \RCA\, suffices to prove that $(P, \leq_Q)$ has property $\sP$.
\end{lem}

\proof See \cite{Marcone20}. \qed

The statement of the previous lemma for {\wqo}(ext) is still open, see \cite[Question 2.15]{Marcone20}; similarly, it is still open in the context of constructive mathematics \cite{BSB23}. In the case of the product the situation is far more complex, as showed by the next result.
\begin{theor}[Cholak, Marcone, Solomon]\label{theor:wqoproduct} Let $\sP$ be any of the properties {\wqo}, {\wqo}(anti), or {\wqo}(ext). The following holds:

\begin{itemize}

\item if $(P,\leq_P)$ and $(Q,\leq_Q)$ are {\wqo}(set), then \RCA\, suffices to prove that $(P \times Q, \leq_{P  \times Q})$ is {\wqo}(set);

\item if $(P,\leq_P)$ and $(Q,\leq_Q)$ have property $\sP$, then \textsf{WKL$_0$} does not suffice to prove that $(P \times Q, \leq_{P  \times Q})$ has property $\sP$. 

\end{itemize}

\end{theor}

\proof See \cite{CMS04,Marcone20}. \qed

These limitations regarding the product of {\wqo} will partially affect the proof of \Theor{theor:Equivw}. For sake of completeness, let us observe that the property ``$P,Q$ {\wqo} implies $P \times Q$ {\wqo}'' is known to be strictly between the principles CAC and ADS \cite[Theorem 1.16]{Towsner20}. We consider the preservation of {\wqo} with respect to functions, starting with the following definition.
\begin{defi}\label{defi:ordermaps} Given two qo $(P,\leq_P)$ and $(Q,\leq_Q)$, a function $\phi \! : P \frec Q$~is:

\begin{enumerate}

\item \emph{order-preserving} if  $p_1 \leq_P p_2 $ implies $  \phi (p_1) \leq_Q \phi (p_2)$;

\item \emph{order-reflecting} if $\phi(p_1) \leq_Q \phi(p_2) $ implies $ p_1 \leq_P p_2$;

\item an \emph{order isomorphism} if it is an order-preserving, order-reflecting bijection.

\end{enumerate}

\end{defi}
Regarding {\wqo} preservation with respect to such order functions, we have the following result whose provability over \RCA\, is, to the best of our knowledge, novel in the literature.
\begin{prop}\label{theor:wqorder} Let $\sP$ be any of the properties {\wqo}, {\wqo}(set), {\wqo}(anti) or {\wqo}(ext), let $P$ and $Q$ be two qo and let $\phi \! : P \frec Q$ be a function; then

\begin{enumerate}

\item if $P$ has the property $\sP$ and $\phi$ is an order-preserving surjection, then $Q$ has the property $\sP$;

\item if $Q$ has the property $\sP$, with $\sP \neq $ {\wqo}(ext), and $\phi$ is an order-reflecting map, then $P$ has the property $\sP$.

\end{enumerate}

\end{prop}

\proof We start with the first point, i.e., $\phi$ order-preserving surjection, considering all the cases for $\sP$.

\smallskip

\noindent \ita{[$P$ {\wqo}]} Let $(q_k)_k$ be an infinite sequence in $Q$, since $\phi$ is surjective we can find (also in \RCA) an infinite sequence $(p_k)_k$ in $P$ such that for all $k$, $\phi(p_k)=q_k$; but $P$ is {\wqo}, thus there exist indexes $i<j$ such that $p_i \leq_P p_j$ and, since $\phi$ is order-preserving, we obtain $q_i=\phi(p_i) \leq_Q \phi(p_j)=q_j$.

\smallskip

\noindent \ita{[$P$ {\wqo}(set) or $P$ {\wqo}(anti)]} The previous technique works here too.

\smallskip

\noindent \ita{[$P$ {\wqo}(ext)]} In order to keep the notation readable, we only consider the case $P,Q$ partial orders, although the following strategy extends to quasi-orders. Let $\preceq_Q$ be a linear extension of $\leq_Q$, we have to prove that $\preceq_Q$ is well-founded. For each element $q \! \in \! Q$, we consider the subset $P_q= \phi^{-1}(q) =\{ p\! \in \! P \, |\, \phi(p)=q \}$ together with the restriction $\leq_P \downharpoonright_{P_q}$ of $\leq_P$ to $P_q$; by Szpilrajn`s theorem, available in \RCA (\cite[Observation 6.1]{DHLS03}, see also \cite{FM12}), we can consider for each $q\! \in \! Q$ 
a linear extension $\tri_q$ of $(P_q, \leq_P \downharpoonright_{P_q})$. We define the following order on $P$:
\[
p_1 \preceq_P p_2 \defAs \phi(p_1) \prec_Q \phi(p_2) \, \vel \, (\phi(p_1)=\phi(p_2) \et p_1 \tri_{\phi(p_1)} p_2).
\]
 It is straightforward to verify that $\preceq_P$ is a linear extension of  $\leq_P$ and $p_1 \! \preceq_P \! p_2$ implies $\phi(p_1) \! \preceq_Q \! \phi(p_2)$.

Since $\preceq_P$ is a linear order extending $\leq_P$ and $P$ is {\wqo}(ext), then $\preceq_P$ is a well-order. Let us assume, by contradiction, that $q_1 \succ_Q q_2 \succ_Q \dots $ is a strictly descending chain in $(Q, \preceq_Q)$, since $\phi$ is surjective we can find a sequence $p_1, p_2, \dots$ such that $\phi(p_k)=q_k$ for all $k$; but $\preceq_P$ is a well-order, thus there exist $i \! < \! j$ such that $p_i \preceq p_j$ and then $\phi(p_i) \preceq_Q \phi(p_j)$, contradiction.

\smallskip

Let us consider the second point, namely $\phi$ order-reflecting.

\noindent \ita{[$Q$ {\wqo}]} Let $(p_k)_k$ be an infinite sequence in $P$, then $(\phi(p_k))_k$ is an infinite sequence in $Q$ and thus there exist indexes $i\! < \! j$ such that $\phi(p_i) \! \leq_Q \! \phi(p_j)$; but, since $\phi$ is order-reflecting, this implies $p_i \! \leq_P \! p_j$.

\smallskip

\noindent \ita{[$P$ {\wqo}(set) or $P$ {\wqo}(anti)]} The previous strategy works here too. \qed

\noindent As before, the case {\wqo}(ext) is still open for an order-reflecting map $\phi$. 

\smallskip

We conclude this section by formally stating Higman's lemma, recalling that $Q^* \putAs \bigcup_{n=0}^{\infty}Q^n$ is the set of finite sequences on $Q$.
\begin{lem}	\label{lem:Higman} If $(Q,\leq)$ is a {\wqo}, then $(Q^*,\leq^*)$ is a {\wqo}, where:
\[
p_1\dots p_n \leq^* q_1 \dots q_m \ \ \mbox{iff} \ \ \exists 1 \! \leq \! i_1 \! < \! \dots \! < \! i_n \! \leq \! m\ \forall k\! \in \! \{1,\dots,n\}\ p_k \leq q_{i_k} .
\]
\end{lem}
Higman's lemma, together with Dickson's lemma, will be analysed from a proof-theoretic point of view in \Sec{sec:HD}; for now, let us recall that, provably over \RCA, Higman's lemma is equivalent to the formal system \ACA.

Remarkably, if in Higman's lemma finite sequences are substituted with infinite sequences (and an appropriate embedding relation), then the corresponding statement is false, as proved by Rado's counterexample \cite{Rado54}. The fact that not all operations preserve the property of being a {\wqo}, particularly the infinitary ones, is the main reason that motivated the introduction of \ita{better quasi-orders} \cite{CP20,Nash-W68}.

\section{Abstract Algebras and Trees}\label{sec:algebras}

In this section, we introduce the two main topics of the paper, Higman's and Kruskal's theorem, as well as the required preliminaries.

\subsection{Ordered Algebras and Higman's Theorem}\label{subsec:algebras}

This paragraph is devoted to formally present \ita{ordered algebras} and \ita{Higman's theorem}. As before all the definitions are given in \RCA; however, since dealing in \RCA\, with sets of operations, i.e., sets of functions and thus sets of sets, required some care, the following definition is slightly more involved than its classical counterpart \cite{Higman52}.
\begin{defi}\label{defi:n-algebra}
An $n$-\emph{ary (abstract) algebra} $(A,M)$ is given by a set $|A| \! \subseteq \! \mN$ together with a finite list $M=(M_i)_i$ with $1\! \leq \! i \! \leq \! n$ of sequences, possibly finite, infinite or empty, of operations over $|A|$; for each $i$, $M_i$ is the sequence of $i$-ary operations from $|A|^i$ to $|A|$.
\end{defi}
As for {\qo}, we tacitly write $A$ instead of $|A|$, omitting also the term ``abstract''. Moreover, although it does not formally exist in \RCA\, as a set, we use $M$ to refer to the family of all the operations of the algebra; namely $\mu \in M$ is not a well-defined formula of our syntax, but just a syntactic abbreviation of ``$\mu$ is an operation over $A$''.

\begin{rem}\label{rem:codes} Since each operation $\mu$ is uniquely determined by its arity $r$, i.e., $\mu \! \in \! M_r$ for some $r\! \in \! \{1, \dots, n\}$, and its position $d$ in the sequence $M_r$, we can uniquely assign a natural number $|\mu|$ to each operation, allowing us to consider also the set $|M_i|$ of the codes of elements of $M_i$; moreover, we can also require that codes for operations and codes for elements of $A$ are disjoint. Given this bijection, at least in the metatheory, in the following we write $\mu$ for both an operation and the corresponding code relying on the context to make clear the exact meaning; similarly for $M_i$ and $|M_i|$. 
\end{rem}
The following generalization can be stated:
\begin{defi}\label{defi:algebra}
An \emph{(abstract) algebra} $(A,M)$ is given by a set $|A| \! \subseteq \! \mN$ together with a sequence $M$ of operations over $|A|$, to any operation $\mu \in M$ is associated its arity $d(\mu)$.
\end{defi}
As before, we simply write $A$ instead of $|A|$, omitting the word ``abstract'' and, frequently, also the arity operator $d$. In analogy with the previous case, we can associate a unique code to each operation of $M$ (and thus considering the code set $|M|$), requiring these codes to be different from the codes of elements of $A$. Finally, we underline how an $n$-ary algebra can be seen as an algebra placing $M=M_1 \cup \dots \cup M_n$ and $d(\mu)=i$ if $\mu \! \in \! M_i$.

For algebras, there are the following auxiliary definitions.
\begin{defi}\label{defi:subalgebr}
Given an algebra $(A,M)$:

\begin{itemize}

\item $(B,M)$ is a \emph{subalgebra} of $(A,M)$ if $B \subseteq A$ and $B$ is closed with respect to $M$, i.e., for every $\mu$ in $M$ of arity $d(\mu)$ and every $b_1, \dots , b_{d(\mu)}$ in $B$, $\mu(b_1, \dots, b_{d(\mu)})  \in B$;

\item $C \subseteq A$ is a \emph{generating set} of $A$, and $A$ is \emph{generated} by $C$, if for every subalgebra $(B,M)$ of $(A,M)$, $C \subseteq B$ implies $B=A$.

\end{itemize}

\end{defi}
Quasi-orders and abstract algebras can now be linked together.

\begin{defi}\label{def:OrderedAlgebra} Given an algebra $(A,M)$ and a quasi-order $\leq$ over $A$:

\begin{enumerate}

\item $(A,M,\leq)$ is an \emph{ordered algebra} if for any function $\mu \! \in \! M$ with arity $d(\mu)$
\[
\forall i\! \leq \! d(\mu) \ a_i \! \leq \!  b_i  \ \ \mbox{implies}\ \ \mu(a_1, \dots, a_{d(\mu)}) \! \leq \! \mu(b_1, \dots, b_{d(\mu)});
\]

\item $\leq$ is a \emph{divisibility order} if, in addition, for any function $\mu \! \in \! M$ with arity $d(\mu)$ 
\[
\forall i\! \leq \! d(\mu) \ \, a_i \leq \mu(a_1, \dots, a_{d(\mu)});
\]

\item $\leq$ is \emph{compatible} with a qo $\tri$ on $M$ if for any $\mu, \la$ in $M$
\[
a_1\dots a_{d(\la)} \! \leq^* \! b_1\dots b_{d(\mu)}\ \mbox{and}\ \la \! \tri \! \mu \ \, \mbox{imply}\ \, \la(a_1, \dots, a_{d(\la)}) \! \leq \! \mu(b_1,\dots , b_{d(\mu)}), 
\]
where $\leq^{*}$ is Higman's order (see \Lem{lem:Higman}).
\end{enumerate}

\end{defi}
This definition can be given almost untouched for an $n$-ary algebra too, the only concrete modification is for compatibility. In an $n$-ary algebra, we could have different qo's $\tri_r$, one for each operation set $M_r$; compatibility is then stated separately for each arity, i.e., if $\la, \mu$ belong to $M_r$ then
\[
\forall i\! \leq \! r \ a_i \! \leq \!  b_i \ \ \mbox{and}\ \ \la \! \tri_r \! \mu \ \ \mbox{imply}\ \ \la(a_1, \dots, a_r) \! \leq \! \mu( b_1, \dots, b_r).
\]
Properties \ita{1.} and \ita{2.} of the previous definition, or their equivalents, can be found also in other contexts. A function over ordinals satisfying those conditions is commonly called, respectively, \ita{monotonic} and \ita{increasing} \cite[def. 4.5]{Schmidt79}; while, talking about natural numbers, such a function is dubbed, respectively, \ita{weakly increasing} and \ita{inflationary} \cite[p.178]{Weiermann06PA}.

\begin{rem}
(cf. \Rem{rem:codes}) Formally speaking in the third point of \Def{def:OrderedAlgebra} (compatibility),  given a qo $\tri$ over $|M|$, we are requiring that $\leq$ preserves this order, namely
\[
a_1\dots a_{d(\la)} \! \leq^* \! b_1\dots b_{d(\mu)}\ \mbox{and}\ |\la| \! \tri \! |\mu| \ \mbox{imply}\ \la(a_1, \dots, a_{d(\la)}) \! \leq \! \mu(b_1,\dots , b_{d(\mu)}).
\]
\end{rem}
The finite and the infinite versions of Higman's theorem can finally be stated; the former is found in the seminal paper of Higman \cite{Higman52}, while the latter is presented in \cite{Pouzet85}. For each $n \! \in \! \mN$, we have the corresponding version of the following theorem:

\begin{theor}[$n$-finite Higman's theorem]\label{theor:finiteHigman} Given an $n$-ary ordered algebra $(A,M,\leq)$, if $(M_r, \tri_r)$ is a {\wqo} for each $r$, $\leq$ is a divisibility order compatible with the qo's of $M$ and $(A,M)$ is generated by a {\wqo} set $C\subseteq A$, then $(A,\leq)$ is {\wqo}.

\end{theor}
Using a notation similar to the one reserved for the various versions of Kruskal's theorem (see below), each $n$-version of the previous theorem is denoted by HT$(n)$, e.g. HT$(3)$ corresponds to the theorem for $3$-ary ordered algebras. Given this notation, the original theorem due to Higman \cite{Higman52} simply reads as follows.

\begin{theor}[Original Higman's theorem]\label{theor:originalHigman} $\forall n\, \mbox{HT}(n)$.
\end{theor}
Using abstract algebras, an infinite version of Higman's theorem, denoted by HT$(\omega)$ and firstly proposed by Pouzet \cite[Theorem 2.2]{Pouzet85}\footnote{Who traces it back to Cohn \cite{Cohn81}.}, can be stated.

\begin{theor}[Infinite Higman's theorem]\label{theor:infiniteHigman} Given an ordered algebra $(A,M,\leq)$, if $(M, \tri)$ is a {\wqo}, $\leq$ is a divisibility order compatible with $\tri$ and $(A,M)$ is generated by a {\wqo} set $C\subseteq A$, then $(A,\leq)$ is {\wqo}.
\end{theor}
The crucial difference between the finite and the infinite version is the presence in the former of an arbitrarily high, but fixed, upper bound on the arity of the operations in $M$. The main references for abstract algebras are Higman's original paper \cite{Higman52} and the book of Cohn \cite{Cohn81}.

\subsection{Labelled Trees and Kruskal's Theorem}\label{subsec:trees}

We start by recalling some definitions regarding trees, tree embedding, and branching degree; as before, all definitions are stated in \RCA.

\begin{defi}\label{def:trees} Given a set $Q$, we inductively define the set $\mT(Q)$ of (finite ordered) \emph{trees} with \emph{labels} in $Q$ as follows:

\begin{enumerate}

\item for each $q \! \in \! Q$, $q[]$ is an element of $\mT(Q)$;

\item if $t_1, \dots , t_n$, with $n \! \geq \! 1$, is a finite sequence of elements of $\mT(Q)$ and $q \! \in \! Q$, then 
$t \putAs q[t_1, \dots , t_n]$ is an element of $\mT(Q)$.

\end{enumerate}

\end{defi}
Since all trees considered in this paper are finite and ordered, we omit these specifications. Connected to labelled trees, there are the following definitions.
\begin{defi}\label{def:labels} Let $Q$ and $\mT(Q)$ be as before:

\begin{itemize}

\item if $t = q[t_1, \dots , t_n]$, then $l(t) \putAs q$ is the \emph{label} of the $t$ and $t_1, \dots , t_n$ are called the \emph{(immediate) subtrees} of $t$;

\item the set of \textsc{nodes} $\sN(t)$ of $t \in \mT(Q)$ is recursively defined as:

       \begin{enumerate}
       
       \item if $t=q[]$, then $\sN(t) \putAs \{ t\}$,
       
       \item if $t=q[t_1, \dots, t_n]$, then $\sN(t) \putAs \{ t\} \cup \sN(t_1) \cup \dots \cup \sN(t_n)$;
       
       \end{enumerate}

\item if $Q$ is a singleton, then $\mT \putAs \mT(Q)$ is the set of \emph{unlabelled} trees.

\end{itemize}

\end{defi} 
Starting from a qo $Q$, it is possible to define an embeddability relation.
\begin{defi}\label{def:treemb} Given a qo $(Q,\leq_Q)$ and $t,s \! \in \! \mT(Q)$, we inductively define the embeddability relation $\preceq$ on $\mT(Q)$ as follows: $t \preceq s$ holds if

\begin{enumerate}

\item $t=p[]$, $s=q[]$ and $ p \leq_Q q$; or

\item $s=q[s_1, \dots , s_m]$ and $t \preceq s_i$ for some $ 1 \! \leq \! i \! \leq \! m$; or

\item $t=p[t_1, \dots, t_n]$, $s=q[s_1, \dots , s_m]$, $p \leq_Q q$ and $t_1\dots t_n \preceq^{*} s_1 \dots s_m$.

\end{enumerate}
Here, $\preceq^{*}$ is Higman's order with respect to $\preceq$, see \Lem{lem:Higman}.
\end{defi}
Let us observe that $(\mT(Q), \preceq)$ is a quasi-order.

Before stating Kruskal's theorems, one last definition is required: the branching degree.

\begin{defi} The \textsc{degree}, $d$, and the \textsc{branching degree}, $Deg$, of a tree $t \! \in \! \mT(Q)$ are recursively defined as:
\begin{enumerate}

\item if $t=q[]$, then $d(t)=Deg(t) \putAs 0$;

\item if $t=q[t_1, \dots, t_n]$, then $d(t) \putAs n$ and $Deg(t) \putAs \max\{n, Deg(t_1), \dots, Deg(t_n)\}$.

\end{enumerate}
Moreover, for every $n \! \in \! \mN$ and every set $Q$:

\begin{itemize}

\item $\mT_n$ is the set of unlabelled trees with branching degree less or equal to~$n$;

\item $\mT_n(Q)$ is the set of trees with labels in $Q$ and branching degree less or equal to $n$.

\end{itemize}

\end{defi}
Since only finite trees are considered, $d$ and $Deg$ are always well defined.

We now present the various versions of Kruskal's theorem under scrutiny; each is denoted by an abbreviation, possibly accompanied by a numeric parameter.

\begin{theor}\label{theor:Kruskal} The following statements hold for every natural number $n$:

\begin{description}

\item[\KTn:] $(\mT_n, \preceq)$ is a {\wqo};

\item[\KTw:] $(\mT,\preceq)$ is a {\wqo};

\item[\KTln:] if $Q$ is a {\wqo}, then $(\mT_n(Q), \preceq)$ is a {\wqo};

\item[\KTlw:] if $Q$ is a {\wqo}, then $(\mT(Q), \preceq)$ is a {\wqo}.

\end{description}

\end{theor}
Starting from the original paper by Kruskal \cite{Kruskal60}, proofs of Kruskal's theorem, or its variation, are ubiquitous in the literature \cite{Goubault13,Nash-W63,Sternagel14,Veldman04}. For what concerns its strength from the reverse mathematics point of view, by a result due to Harvey Friedman \cite{Simpson85nonprovability}(see also \cite{Gallier91}), it is known that even Kruskal's theorem for unlabelled trees (\KTw \ in our notation), as well as its finite miniaturization \cite{Simpson85nonprovability,Smith85}, is not provable in \textsf{ATR$_0$}. A thorough investigation of the unlabelled versions of Kruskal's theorem has been carried out by Rathjen and Weiermann \cite{RW93}, some of their results can be summarized in the following theorem.

\begin{theor}[Rathjen and Weiermann]\label{theor:Kruskal}
\[
\mbox{RCA}_0 \vdash \, \forall n\, \mbox{KT}(n) \, \sse \, \mbox{KT}(\omega) \, \sse \,  \mbox{WO}\left( \vartheta (\Omega^{\omega}) \right).
\]
\end{theor}
In the previous statement, WO$\left(\vartheta ( \Omega^{\omega}) \right)$ denotes the well-orderedness of the countable ordinal $\vartheta (\Omega^{\omega})$, where $\vartheta$ is a so-called \ita{collapsing function}, see \cite{Buchholz86,RW93,VanderMereen15} for an introduction to such functions. As we will see, the situation for the labelled case is quite different.

\section{Relations between Higman's and Kruskal's Theorems}\label{sec:HK}

Adapting to the framework of reverse mathematics the proof by Pouzet \cite{Pouzet85} (see also \cite{Schmidt79}), we prove that, regarding ordered algebras and labelled trees, Higman's and Kruskal's theorems are two sides of the same coin. We first consider the fine case, namely HT$(n)$ and KT$_\ell(n)$; our goal is to prove the following equivalence.
\begin{theor}\label{theor:equivn} RCA$_0 \ \vdash \ \forall n\, HT(n) \, \longleftrightarrow \,  \forall n\, KT_\ell(n)$.
\end{theor}
The previous result is an immediate consequence of the following lemma that is provable in \RCA.
\begin{lem}\label{lem:foralln}\label{foralln} $\forall n$ HT$(n) \! \longleftrightarrow \!$ KT$_\ell(n)$.
\end{lem}

\proof Let us prove \Lem{lem:foralln} by fixing $n$ and considering each direction separately. First of all, we introduce some notation:
\begin{itemize}

\item $ \leq$ indicates the {\qo} over the set $A$, representing an $n$-ary algebra, whose elements are denoted by $ a,b,c \dots $; $\leq$ will be also the usual order over~$\mN$.

\item $\tri$ stands for the {\qo} over $M_r$ for each $r$, whose elements are denoted by $\mu, \la, \dots$; for sake of readability, we use the same symbol $\tri$ for each set $M_r$, $1 \! \leq \! r \! \leq \! n$.

\item $\preceq$ represents the embeddability relation between finite ordered labelled trees, which are denoted by $t,s, \dots $; if necessary, also the extended notation $t=q[t_1, \dots , t_n]$ is used.

\end{itemize}

\noindent ``HT$(n) \imp \mbox{KT}_\ell(n)$''

\smallskip

Given a {\wqo} $(Q,\leq)$, we must prove that $(\mT_n(Q), \preceq)$ is well quasi-ordered by~$\preceq$.
For each $1 \leq  r  \leq  n$ and $q \! \in \! Q$, we define the following operation over $\mT_n(Q)$:
\[
\bigoplus_{r,q} : \mT_n(Q)^r \frec \mT_n(Q)\ \mbox{with} \ \bigoplus_{r,q}(t_1, \dots, t_r) \putAs q[t_1, \dots , t_r]
\]
where $q[t_1, \dots , t_r]$ is the tree whose root is labelled by $q$ and whose immediate subtrees, are $t_1, \dots, t_r$.
We dub $\bigoplus_{N,Q}$ the family of all these operations (which is in bijection with $\{1, \dots, n\} \times Q$); for each of these operations we can consider a unique code (see \Rem{rem:codes}). For each $1 \leq  r  \leq  n$, we quasi-order the set $\bigoplus_{r,Q}=\{ \bigoplus_{r,q}\, | \, q \! \in \! Q \}$ (of codes) of $r$-ary operations in the following way: $\bigoplus_{r,q} \tri_r \bigoplus_{r,p}$ iff $q \! \leq \! p$; with such quasi-orders each $(\bigoplus_{r,Q}, \tri_r)$ is a {\wqo}, since also $(Q, \leq)$ is a {\wqo}.

\smallskip

\noindent \ita{CLAIM.} $(\mT_n(Q), \bigoplus_{N,Q}, \preceq)$ is an ordered $n$-ary algebra with $\preceq$ a divisibility order compatible with $\tri$; moreover $C=\{q[]\, |\, q \! \in \! Q\}$ is a generating set for~$\mT_n(Q)$.

\noindent \ita{Proof of the Claim}: Since $\bigoplus_{N,Q}$ is a family of at most $n$-ary operations, we need to check the three properties of \Def{def:OrderedAlgebra} regarding the quasi-order $\preceq$ and the generation of $\mT_n(Q)$ by $C$.

\smallskip

\noindent \emph{Ordered algebra:} given $1 \! \leq \! r \! \leq \! n,\ q \! \in \! Q$ and $t_i \! \preceq \! s_i$ for all $1 \! \leq \! i \! \leq \! r$, then $q[t_1, \dots, t_r] \preceq q[s_1, \dots, s_r]$ holds by the third condition of \Def{def:treemb}; thus, $\bigoplus_{r,q}(t_1, \dots , t_r) \! \preceq \! \bigoplus_{r,q}(s_1, \dots , s_r)$.

\smallskip

\noindent \emph{Divisibility:} given $1 \! \leq \! r \! \leq \! n$ and $q \! \in \! Q$, $t_i \! \preceq \! q[t_1, \dots, t_r]$ holds by the second condition of \Def{def:treemb}; therefore, $t_i \! \preceq \! \bigoplus_{r,q}(t_1, \dots , t_r)$ for all $1 \! \leq \! i \! \leq \! r$.

\smallskip

\noindent \emph{Compatibility:} given $1 \! \leq \! r \! \leq \! n$, $\bigoplus_{r,q} \! \tri \! \bigoplus_{r,p}$ and $t_i \!  \preceq  \! s_i$ for all $1 \! \leq \! i \! \leq \! r$, then $q \! \leq \! p$ and $q[t_1, \dots, t_r] \! \preceq \! p[s_1, \dots, s_r]$ holds by the third case of \Def{def:treemb}; thus $\bigoplus_{r,q}(t_1, \dots, t_r) \preceq \bigoplus_{r,p}(s_1, \dots, s_r)$.

\smallskip

\noindent \emph{Generation:} obviously $C\! \subseteq \! \mT_n(Q)$; we need to prove that, given a subset $B\! \subseteq \! \mT_n(Q)$ closed with respect to $\bigoplus_{N \!,Q}$, if $C \! \subseteq \! B$, then $B = \mT_n(Q)$. Let $C \! \subseteq \! B$ and $t \! \in \! \mT_n(Q)$, if $t \! \in \! C$, i.e., $t=q[]$ with $q \! \in \! Q$, then $t \! \in \! B$; else $t$ is obtained by a finite number of elements of $C$, the leaves of $t$, and a finite number of applications of operations in $\bigoplus_{N,Q}$. Since $B$ is closed with respect to $\bigoplus_{N,Q}$, $t \! \in \! B$ and thus $B = \mT_n(Q)$. Lastly, $(C, \preceq )$ is {\wqo} since $q[] \! \preceq \! p[] \Leftrightarrow q \leq p$ and $(Q,\leq)$ is {\wqo} by hypotheses. \ita{End of the proof of the Claim}.

\smallskip

To conclude, since $\tri$ is a {\wqo} on $\bigoplus_{N,Q}$, we can apply HT$(n)$ to obtain that $(\mT_n(Q), \preceq)$ is a {\wqo}.

\bigskip

\noindent ``KT$_\ell(n) \imp \mbox{HT}(n)$''

\smallskip

Let $(A,M,\leq)$ be an $n$-ary ordered algebra satisfying all the hypotheses of Higman's theorem: each set of $r$-ary operations $M_r$ is a {\wqo} with respect to $\tri$; $\leq$ is a divisibility order compatible with the {\wqo} $\tri$ of each $M_r$; there is a {\wqo} generating set $C$. We must prove that $(A,\leq)$ {\wqo}.

Since $M$ is a {\wqo} (given by the finite disjoint union $M_1 \cup \dots \cup M_n$ of {\wqo}'s) and $C$ is a {\wqo}, so is $C\! \cup \! M$, by \Lem{lem:wqounion}, and applying \KTln \ so is $(\mT_n(C\! \cup \! M), \preceq)$. Over $C\! \cup \! M$, we consider the \ita{grade function} $g\! :\! C \!  \cup \! M \frec \mN$ defined as $g(x) \putAs 0$ if $x \! \in \! C$ and $g(x) \putAs r$ if $x\! \in \! M_r$. Finally,  we define, inside $\mT_n(C\! \cup \! M)$, the subset $\mT_n^g(C\! \cup \! M)$ of \ita{graded trees}, i.e., the trees $t$ satisfying the following grade condition
\[
\forall s\! \in \! \sN(t) \ \ d(s)=g(l(s)).
\]
In other words, we are considering trees $t$ such that, for each node $s$, the degree $d(s)$ of the node is equal to the grade $g$ of the label $l(s)$ of the node itself. The functions $g,d$ and $\mT_n^g(C\! \cup \! M)$ are all definable in \RCA; moreover, $\mT_n^g(C\! \cup \! M)$ is well quasi-ordered by $\preceq$.

We now construct an order-preserving function $\phi$ between $(\mT_n^g(C\! \cup \! M), \preceq)$ and $(A,M,\leq)$ that is defined by recursion over the tree structure of the elements of $\mT_n^g(C\! \cup \! M)$:
\[
\phi (t) = \left\{
  \begin{array}{ll}
     q & \mbox{if}\ t=q[]\ \mbox{with}\ q\! \in \! C;\\
     \mu (\phi(t_1), \dots , \phi(t_r)) & \mbox{if}\ t=\mu[t_1, \dots , t_r]\ \mbox{with}\ \mu \! \in \! M_r.
   \end{array}  
   \right.
\]
Since we are considering graded trees, the label of a leaf is always an element of $C$; whereas the label of a node of degree $r$ is always (the code of) an operation of $(A,M)$ of arity $r$, i.e., an element of $M_r$. This, together with the fact that we are treating ordered trees (i.e., $\la[t_1,t_2] \not= \la[t_2,t_1]$, unless $t_1=t_2$), ensures that $\phi$ is well defined. An instance of $\phi$ is depicted below.

\begin{equation*}
  \begin{tikzpicture}[baseline={([yshift=-.5ex]current bounding box.center)}, grow'=up, level distance=10mm, sibling distance=10mm, vertex/.style={anchor=base, circle,fill=black!25,minimum size=18pt,inner sep=2pt}]
    \usetikzlibrary {quotes}
       \node {$\mu$ } 
        child {node {a}}
        child {node {$\la$}
        child {node {b}}
        };
  \end{tikzpicture}
  \ \ \longmapsto \ \ \mu(a, \la(b))
\end{equation*}
We prove by induction that $\phi$ is order-preserving, i.e., $t\! \preceq \! s \imp \phi(t) \!\leq \! \phi(s)$, considering the three defining cases of \Def{def:treemb}. Given $t,s \! \in \! \mT_n^g(C\! \cup \! M)$:
\begin{enumerate}

\item if $t \! \preceq \! s$ with $t=q[], s=p[]$, then $q \leq p$ and thus $\phi(t) \! \leq \! \phi(s)$; \vspace{0.1cm}

\item if $t \! \preceq \! s$ with $s=\mu[s_1, \dots , s_r]$, $\mu \! \in \! M_r$ and $t \preceq s_i$ for some $i\! \in \! \{ 1, \dots, r \}$, then, by induction, $\phi(t) \! \leq \! \phi(s_i)$ and, since $\leq$ is a divisibility order, $\phi(t) \! \leq \! \mu( \phi(s_1), \dots, \phi(s_r))$, thus $\phi(t) \! \leq \! \phi(s)$;
\vspace{0.1cm}

\item if $t \! \preceq \! s$ with $t=\mu [t_1, \dots , t_r]$, $s=\la [s_1, \dots , s_k]$, $\mu \! \in \! M_r$, $\la \! \in \! M_k$, $\mu \tri \la$ and there exist $1 \! \leq \! i_1 \! < \! \dots \! < \! i_r \! \leq \! k$ such that $t_1 \! \preceq s_{i_1}, \dots , t_r \! \preceq \! s_{i_r}$, then, by definition of $\tri$, $r=k$ and, by induction, $\phi(t_1) \! \leq \! \phi(s_1), \dots , \phi(t_r) \! \leq \! \phi(s_r)$, thus, since $\leq$ is compatible with $\tri$, $\mu(\phi(t_1), \dots ,\phi(t_r)) \! \leq \! \la(\phi(s_1), \dots, \phi(s_r))$, i.e., $\phi(t) \! \leq \! \phi(s)$.

\end{enumerate}

In order to conclude by applying \Theor{theor:wqorder}, we need to prove that $\phi$ is also surjective. Let us consider $\phi(\mT_n^g(C\! \cup \! M)) \! \subseteq \! A$. Obviously $C \subseteq \phi(\mT_n^g(C \cup  M))$; since $C$ is a generating set it suffices to check that $\phi(\mT_n^g(C \cup  M))$ is closed under $M$. Given $a_1, \dots, a_r \! \in \! \phi(\mT_n^g(C\! \cup \! M))$ and $\mu \! \in \! M_r$, there exist $t_1, \dots  , t_r \! \in \! \mT_n^g(C\! \cup \! M)$ such that $\phi(t_1)=a_1, \dots , \phi(t_r)=a_r$; but now $\mu[t_1, \dots, t_r]$, i.e., the tree with the root labelled by $\mu$ and with immediate subtrees $t_1, \dots , t_r$, is an element of $\mT_n^g(C\! \cup \! M)$ and, by definition, $\phi(\mu[t_1, \dots, t_r])=\mu(\phi(t_1), \dots , \phi(t_r))=\mu(a_1, \dots, a_r)$, thus $\mu(a_1, \dots , a_r) \! \in \! \phi(\mT_n^g(C\! \cup \! M))$. So $\phi(\mT_n^g(C\! \cup \! M))=A$, i.e., $\phi$ is surjective, and we can apply \Theor{theor:wqorder} to conclude that $(A,\leq)$ is {\wqo}. \qed

\smallskip

We move to the infinite case of the previous equivalence.

\begin{theor}\label{theor:Equivw} \RCA\, $\vdash$ HT($\omega$) $\longleftrightarrow$ KT$_\ell$($\omega$).
\end{theor}

\proof The general framework is quite similar, the main differences regarding the background are the following:
\begin{itemize}

\item We are dealing with abstract algebras rather than $n$-ary algebras; this means that we could have operations of any finite arity. In particular, since an infinite union of {\wqo}'s is not in general a {\wqo}, we do not consider a {\wqo} $\tri_r$ for each set $M_r$ of $r$-ary operations as before; conversely, we have one single {\wqo} $\tri$ over the set $M$.

\item We are dealing with trees of, possibly, any finite branching degree, i.e., we consider $\mT(Q)$ and $\mT^g(Q)$ rather than $\mT_n(Q)$ and $\mT_n^g(Q)$.

\end{itemize}
For what concerns the proof, we briefly highlight the required adaptations being the general strategy essentially the same.

\smallskip

\noindent ``HT$(\omega) \imp \mbox{KT}_\ell(\omega)$'': the main difference involves the operations $\bigoplus_{r,q}$ which are now defined for each natural number $r\! \in \! \mN$. The set (of codes) $\bigoplus_{\mN, Q}$ is well quasi-ordered by the ``quotient order'' over $Q$, i.e., $\bigoplus_{n,q} \! \tri \! \bigoplus_{m,p}$ iff $q\! \leq \! p$ (we can not take directly $\mN \times Q$ since, by \Theor{theor:wqoproduct}, RCA$_0$ does not prove in general that $\mN \times Q$ is a {\wqo}). We have then a corresponding version of the claim.

\smallskip

\noindent \ita{CLAIM.} $(\mT(Q), \bigoplus_{\mN,Q}, \preceq)$ is an ordered algebra with $\preceq$ a divisibility order compatible with $\tri$; moreover $C=\{q[]\, |\, q \! \in \! Q\}$ is a generating set for $\mT(Q)$.

\smallskip

\noindent The proof of the claim, as well as the rest of this direction, is almost equal; a minor modification occurs for compatibility since, through $\tri$, we can compare operations $\bigoplus_{r,q}$ and $\bigoplus_{s,p}$ with different arities, i.e., $r\not= s$.

\bigskip

\noindent ``KT$_\ell(\omega) \imp \mbox{HT}(\omega)$'': in this case, $M$ is directly a {\wqo} and we apply \KTlw \ to obtain the well quasi-orderedness of $(\mT(C \cup M), \preceq)$. The grade function $g$ is now an extension of the arity function $d$ of the abstract algebra under consideration, namely $g\! :\! C \! \cup \! M \frec \mN$ is defined as $g(x) \putAs 0$ if $x \! \in \! C$ and $g(x) \putAs d(x)$ if $x\! \in \! M$. The remaining minor changes concern notation. \qed

\smallskip

From previous equivalences, in particular the one between finite versions, the following corollary is derivable.

\begin{cor} \RCA \,$\nvdash $ KT$(\omega) \, \rightarrow \, \forall n\,$KT$_{\ell}(n)$. 

\end{cor}
\proof Suppose, by contradiction, that the implication is provable over \RCA, then, by \Theor{theor:equivn}, we obtain \RCA\, $\vdash$ KT($\omega$) $ \frec \forall n\,$HT($n$). From $\forall n\,$HT($n$), one can easily deduce Higman's lemma, HL, which correspond to a particular case of HT(2), and so \RCA\, $\vdash$ KT($\omega$) $ \frec$ HL; this means that every model of \RCA\, satisfying KT($\omega$) must satisfy also Higman's lemma. But now, there is a counterexample because REC, the minimal $\omega$-model of \RCA\, given by recursive sets, satisfies KT($\omega$) (since KT($\omega$) is a true $\Pi^1_1$ statement); whereas REC does not satisfy Higman's lemma, for it is not a model of \ACA\, which is equivalent, over \RCA, to Higman's lemma \cite[Theorem X.3.22]{Simpson09}.\footnote{The authors thank Anton Freund for suggesting this approach.} \qed
 
\smallskip 
 
In \cite{BW:ordinals}, the authors computed the $\Pi^1_1$ proof-theoretic ordinals\footnote{See \cite{Rathjen99} for an accessible introduction to ordinal analysis.} of $\forall n \,$\KTln \, and \KTlw \, obtaining the following estimations:

\begin{theor}\label{theor:KTlordinals}
\[
|\RCA + \forall n\, \mbox{KT}_{\ell}(n)|=\vartheta (\Omega^{\omega}+\omega)\ \mbox{and}\ |\textsf{RCA}_0 + \mbox{KT}_{\ell}(\omega)|=\vartheta (\Omega^{\omega+1}).
\] 
\end{theor}
 From this last theorem, the following separation result, which is in striking contrast with \Theor{theor:Kruskal}, derives.

\begin{prop} RCA$_0 \ \nvdash \, \forall n\, \mbox{KT}_{\ell}(n) \ \rightarrow \ \mbox{KT}_{\ell}(\omega)$.
\end{prop}
Previous findings regarding Higman's and Kruskal's theorems can be summed up in the following proof-relation schema where all the possible implications over \RCA\, are depicted.
 
 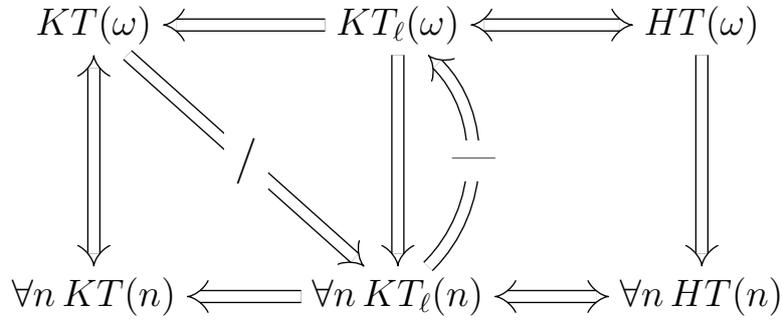
\begin{figure}[!!htb]
 
\begin{center}

\begin{tikzpicture}[every node/.style={rectangle, fill=white, text=black}]
\node (KT)    at (0,1.8)  {\large{$KT(\omega)$}};
\node (KTl)   at (4,1.8)  {\large{$KT_\ell(\omega)$}};
\node (HT)    at (8,1.8)  {\large{$HT(\omega)$}};
\node (KTn)  at (0,-1.8) {\large{$\forall  n\, KT(n)$}};
\node (KTln) at (4,-1.8) {\large{$\forall  n\, KT_\ell(n)$}};
\node (HTn)  at (8,-1.8) {\large{$\forall  n\, HT(n)$}};

\draw[line width=0.5pt, double distance=4pt, -implies] (HT)--(HTn);
\draw[line width=0.5pt, double distance=4pt, -implies] (KTl)--(KTln);

\draw[line width=0.5pt, double distance=4pt,-implies] (KTl)--(KT);
\draw[line width=0.5pt, double distance=4pt,-implies] (KTln)--(KTn);

\draw[line width=0.5pt, double distance=4pt, implies-implies] (KT)--(KTn);

\draw[line width=0.5pt, double distance=4pt, implies-implies] (HT)--(KTl);
\draw[line width=0.5pt, double distance=4pt, implies-implies] (HTn)--(KTln);

\draw[line width=0.5pt, double distance=4pt, -implies] (KT)--(KTln);

\draw[line width=0.5pt, double distance=4pt, -implies] (KTln) edge [double, bend right=45] (KTl);

\node (/)  at (2,0) {\Large{$\slash$}};
\node (//) at (5,0) {\Large{\----}};

\end{tikzpicture}

\end{center}
\caption{\footnotesize Proof-theoretic relations over \RCA\, between finite and infinite versions of Higman's theorem and Kruskal's theorem, with and without labels.}

\end{figure}


\section{Relations between Higman's and Dickson's Lemmas}\label{sec:HD}

Despite having been subsequently overtaken by Kruskal's theorem, Higman's and Dickson's lemma remain two remarkable milestones in the history of well quasi-order theory. In this section, we collect the main proof-theoretic results regarding these two achievements, including their proof-theoretic ordinals. Recall that $Q^{*}$ denotes the set of finite sequences over $Q$.
\begin{defi}\label{defi:HLDL} We adopt the following abbreviations:

\begin{itemize}

\item HL: standard Higman's lemma, i.e., $\forall Q  \, (Q \, {\wqo} \frec Q^* \, {\wqo})$;

\item HL$(\omega)$: $\mN^*$ is {\wqo}, with the order on $\mN$ given by the standard natural order $\leq$.

\item HL$(n)$: $\{0, \dots, n-1\}^*$ is {\wqo}, with the order on $\{0, \dots , n-1 \}$ given by the equality, i.e., $n \! \leq \! m$ iff $n=m$; with side assumption that $n \geq 1 $.

\item DL: Dickson's lemma, i.e., for each $n$, $\mN^n$ is {\wqo} with respect to the product order.
\end{itemize}

\end{defi}
Several findings present in the literature \cite{CMR19,Clote90,Girard87,Hasegawa94,Simpson88basis,Simpson09} are collected in the following  result regarding $\Pi^1_1$ ordinals.
\begin{theor} The following proof-ordinal estimations hold:

\begin{enumerate}

\item $|\RCA + HL |=\ep_0$;

\item $|\RCA +HL(\omega) |=|\RCA + \forall n\, HL(n) |=\omega^{\omega^{\omega+1}} \!$;

\item for all $n$ $|\RCA + HL(n)|= \omega^{\omega^{n}} \!$.

\end{enumerate}

\end{theor}
The proof of the previous estimations derives mainly from the following four lemmas, where WO$(\al)$ denotes the well-ordering of $\al$.
\begin{lem}(Simpson, Girard) \RCA\, $\vdash \mbox{HL} \longleftrightarrow$ \ACA.
\end{lem}
\proof \cite[Theorem X.3.22]{Simpson09}, see also \cite{Girard87,Simpson88basis}. \qed

\begin{lem}(Simpson, Clote) \RCA\, $\vdash \forall n \, \mbox{HL}(n) \longleftrightarrow \mbox{HL}(\omega) \longleftrightarrow \mbox{WO}(\omega^{\omega^{\omega}})$.
\end{lem}
\proof \cite[Theorem 5]{Clote90}, see also \cite{Simpson88basis}. \qed

\begin{lem}(Simpson, Hasegawa) \RCA\, $\vdash \mbox{HL}(n) \longleftrightarrow \mbox{WO}(\omega^{\omega^{n-1}})$.
\end{lem}
\proof In \cite{Hasegawa94} R. Hasegawa proves it improving a well known result due to G. Simpson \cite[sublemma 4.8]{Simpson88basis}. \qed

\begin{lem}(Carlucci, Mainardi, Rathjen) If $\sigma$ is a countable ordinal such that $\omega \x \sigma = \sigma$, then $\left| \RCA + \mbox{WO}(\sigma) \right| = \sigma^{\omega}$. \end{lem}
\proof \cite{CMR19} where, as usual in ordinal analysis, a primitive recursive notation for $\sigma$ is assumed. \qed

\smallskip

Given the above results, we easily obtain.
\begin{cor} Over \RCA\, the following chain of implications is provable:
\[
HL  \frec  HL(\omega)  \leftrightarrow  \forall n\, HL(n) \ \frec \dots \frec HL(n+1) \ \frec \ HL(n)  \frec  \dots \frec HL(2) \leftrightarrow DL.
\]
\end{cor}
The only proof-theoretical relation that does not derive immediately from the aforementioned results is the last one, namely \RCA\, $\vdash \mbox{HL}(2) \leftrightarrow \mbox{DL}$. This can be obtained indirectly through another result by Simpson \cite[Theorem X.3.20]{Simpson09}, i.e., \RCA\, $\vdash \mbox{DL} \leftrightarrow \mbox{WO}(\omega^{\omega})$; nevertheless, we give a direct proof adapting to the framework of reverse mathematics a constructive proof originally proposed by J. Berger \cite{Berger16}. All the following definitions and results are stated or proved in \RCA.
\begin{prop} $HL(2) \leftrightarrow DL$.
\end{prop}
\proof Denoted by $\mN_{+}$ the set of positive integers and by $\{0,1\}^*_0$ the set of finite sequences of $0$ and $1$ starting with a $0$, we actually prove the equivalence over \RCA\, of the two following statements:

\vspace{0.1 cm}

\noindent \textbf{$DL'$:} for all $n$, $\mN_{+}^n$ is a {\wqo};

\vspace{0.1 cm}

\noindent \textbf{$HL'(2)$:} $\{ 0,1\}^*_0$ is a {\wqo};

\vspace{0.1 cm}

\noindent whose equivalence with, respectively, ``$DL$'' and ``$HL(2)$'' is apparent.

\smallskip

We denote by $u,v,w$ the elements of $\{0,1\}_0^*$, by $i$ an element of $\{0,1 \}$, by $n,m,k,l$ the elements of $\mN$ and, with a harmless abuse of notation, by $0$ both the natural number and the string consisting only of the number $0$. Moreover, $ui$ indicates the string obtained from $u$ by adding as final letter $i$, and $|u|$ the length of $u$, i.e., its number of letters. Over $\{0,1\}^{*}_0$ we define the quasi-order $\subeq$ given by subsequences, namely $v \! \subeq \! w$ if $v=v_1 \dots v_p,\ w=w_1\dots w_q$ and there exist $1 \! \leq r_1 \! < \dots <\! r_p \! \leq \! q$ such that $v_j \! = \! w_{r_j}$ for all $1 \! \leq \! j \! \leq \! p$; whereas over $\mN^n_{+}$ let us consider the product order $\leq$, i.e., $(k_1, \dots k_n) \! \leq \! (k'_1, \dots , k'_n)$ if $k_j \! \leq k'_j$ for all $1 \! \leq \! j \! \leq \! n$.

\smallskip

We define the following function $\la \! : \! \{0,1\}_0^* \frec \{0,1\}$ as
\[
  \la(0) \putAs 0 \ \mbox{and} \ \la(ui) \putAs i,
\]
i.e., $\la(w)$ is just the last letter of $w$.

\smallskip

Let us also define $\Phi \! : \! \{0,1\}^{*}_0 \frec \mN$ as
\[
 \Phi(0) \putAs 1 \ \ \mbox{and} \  \Phi(ui) \putAs \left\{
 \begin{array}{ll}
 \Phi(u) & \mbox{if}\ \la(u)=i,\\
 \Phi(u)+1 & \mbox{if}\ \la(u)\not=i.
 \end{array} 
\right.
\]
$\Phi(w)$ is called the \ita{weight} of $w$ and amounts to the number of ``blocks'' of 0's and 1's in $w$, e.g.,
\[
\Phi(00)=1, \ \ \Phi(01100)=3, \ \ \Phi(011011100)=5.
\]

\smallskip

Moreover, we inductively define $F \! :\! \{0,1\}^{*}_0 \frec \bigcup_{n\! \geq \! 1}\mN^n_{+}$ as
\[
F(0) \putAs 1\ \mbox{and, if}\ F(u)=(k_1, \dots, k_m),\ \mbox{then}\ F(ui) \putAs \left\{
 \begin{array}{ll}
 (k_1, \dots, k_m+1) & \mbox{if}\ \la(u)=i,\\
 (k_1, \dots, k_m, 1) & \mbox{if}\ \la(u)\not=i.
 \end{array} 
\right.
\]
$F$ is a bijection which converts a string of weight $m$ into an $m$-tuple counting the occurrences of 0's and 1's, e.g.,
\[
F(00)=(2), \ \ F(01100)=(1,2,2),\ \ F(011011100)=(1,2,1,3,2).
\]
Since we are dealing with finite sequences of natural numbers, $\la, \Phi$, and $F$ are all definable in \RCA.

\smallskip

Fixed $n\! \in \! \mN_{+}$, let $\mathcal{W}_n$ denote the elements of $\{0,1\}^*_0$ with weight $n$ and $F_n$ the restriction of $F$ to $\sW_n$; given the previous definitions we obtain the following claim.

\noindent \ita{CLAIM 1.} $F_n \! :\! \sW_n \frec \mN^n_{+}$ is an order isomorphism between $(\sW_n, \subeq)$ and $(\mN^n_{+},\leq)$.

\noindent \ita{Proof of Claim 1.} It is easy to see that $F_n$ is a bijection since it converts a string of weight $n$ into an $n$-tuple and the other way round. For what concerns the order, we proceed as follows. Let $u \subeq v$ with $u,v \in \mathcal{W}_n$; this means that $u,v$ are made by the same numbers of blocks of $0$'s and $1$'s (starting with $0$) and that both the lengths $|u|$ and $|v|$ are greater or equal to $n$. We adopt the following notation, let $B_j^{\Box}$ and $k_j^{\Box}$ (with $1 \! \leq \! j \! \leq \! n$ and $\Box \! \in \! \{u,v\}$), denote respectively the $j^{th}$-block of $\Box$ and its length (in particular, $B_j^{\Box}$ is a non-empty finite string made entirely either by $0$'s or $1$'s); e.g., $B_3^u$ is the third block of $u$  and $k_3^u$ is the number of $0$'s composing it (since $u$ starts with $0$ and thus the third block, as all the odd ones, is made of $0$'s). Let us observe that $F_n(u)=(k_1^u, \dots, k_n^u)$ and the same for $v$. Since $u \subeq v$ and $u,v \in \mathcal{W}_n$, we obtain $B_j^{u} \subeq B_j^{v}$, and thus $k_j^u \leq k_j^v$, for all $1 \! \leq \! j \! \leq \! n$; namely, $F_n(u)=(k_1^u, \dots, k_n^u) \leq (k_i^v, \dots, k_n^v)=F_n(v)$. For the other direction, if $k=(k_1, \dots, k_n) \leq (k'_1, \dots, k'_n)=k'$, we can consider non-empty finite strings $B_j^k$ and $B_j^{k'}$, with $1 \! \leq \! j \! \leq \! n$, made respectively of $k_j$ and $k'_j$ $0$'s (if $j$ is odd) and $k_j$ and $k'_j$ $1$'s (if $j$ is even). Since $k_j \leq k'_j$, $B_j^k \subeq B_j^{k'}$ and, considering juxtaposition, $F_n^{-1}(k)=B_1^k\dots B_n^k \subeq B_1^{k'}\dots B_n^{k'}=F_n^{-1}(k')$. \ita{This concludes the proof of Claim 1.}

\smallskip

Finally, we define a length-normalizing function $G_n \! :\! \bigcup_{m\! \geq \! 1}\mN^m_{+} \frec \mN^n_{+}$ as
\[
G_n(k_1, \dots, k_m) \putAs \left\{
  \begin{array}{ll}
  (k_1, \dots, k_m, 1, \dots,1) & \mbox{if}\ m<n,\\
  (k_1, \dots, k_n) & \mbox{if}\ n\! \leq \! m.
  \end{array}
  \right.
\]
For $u\! \in \! \{0,1\}^*_0$, let us set
\begin{equation}\label{eq:!_n}
u!_n \putAs F_n^{-1}(G_n(F(u))). 
\end{equation}

\noindent \ita{CLAIM 2.} The following properties hold:
\begin{itemize}

\item $\Phi(u!_n)=n$;

\item if $\Phi(u)=n$, then $u=u!_n$;

\item if $\Phi(u)<n$ (resp. $\Phi(u)>n$), then $|u| < |u!_n|$ (resp. $|u| > |u!_n|$);

\item $G_n(F(u))=F_n(u!_n)$.

\end{itemize}
\noindent \ita{Proof of Claim 2.} Straightforward. 

Before the main result, we state and prove two other claims.

\noindent \ita{CLAIM 3.} For all $v,w$ in $\{0,1\}^*_0$, the following holds
\[
\left( 2 \x |v| \! \leq \! \Phi(w) \right) \ \imp \ v \subeq w
\]

\noindent \ita{Proof of CLAIM 3.} Let us assume, by contradiction, that $v=v_1\dots v_n$, $w=w_1\dots w_m$ with $\Phi(w)=k$, $2n \! \leq \! k$, but $v \! \not  \subeq \! w$. We consider the word $u=01\dots 01$ of length and weight $2n$; since for all $1 \! \leq \! j \! \leq \! n$, $v_j=u_{2j-1}$ or $v_j=u_{2j}$, then $v \! \subeq \! u$. Moreover, since $\Phi(u)=2n\! \leq \! k=\Phi(w)$, we have that $u_j =\bar{w}_j$ where $\bar{w}_j$ is the first letter of the $j$-th block of $w$, and so $u\! \subeq \! w$. Thus $v\! \subeq \! u \! \subeq \! w$, contradiction. \ita{This concludes the proof of Claim 3.}

\noindent \ita{CLAIM 4.} For all $v,w$ in $\{0,1\}^*_0$, we have
\[
\Phi(v) \! \leq \! \Phi(w) \! \leq \! n \, \et \, v!_n \! \subeq \! w!_n \ \imp \ v\! \subeq \! w.
\]

\noindent \ita{Proof of Claim 4.} Let $k=\Phi(v)$. By $v!_n \! \subeq \! w!_n$, we obtain $G_n(F(v)) \! \leq \! G_n(F(w))$ and therefore
\[
F_k(v)=G_k(F(v)) \leq G_k(F(w))=F_k(w!_k),
\]
thus
\[
v\! \subeq \! w!_k.
\]
Since $k\! \leq \! \Phi(w)$, $w!_k \! \subeq \! w$ holds and, by transitivity, $v\! \subeq \! w$. \ita{This concludes the proof of Claim 4.}

\smallskip

We can now prove the main result of this section considering separately the two directions.

\smallskip

\noindent ``$HL'(2) \imp DL'$'': fix $n \! \in \! \mN_{+}$ and an infinite sequence $f\! :\! \mN \frec \mN^n_{+}$. By applying $HL'(2)$ to $F_n^{-1}\circ f$ we obtain $k\! < \! l$ with
\[
F_n^{-1} \circ f(k) \subeq F_n^{-1} \circ f(l),
\]
and thus, applying $F$, $f(k) \! \leq \! f(l)$.

\smallskip

\noindent ``$DL' \imp HL'(2)$'': fix an infinite sequence $g\! :\! \mN \frec \{0,1\}^{*}_0$ and set $n=2 \x |g(0)|$. Then, applying $DL'$, we obtain $k\! < \! l$ such that
\[
(F_n(g(k)!_n),\Phi(g(k))) \leq (F_n(g(l)!_n),\Phi(g(l))).
\]
If $n\! < \! \Phi(g(l))$, then $g(0) \! \subeq \! g(l)$ holds by CLAIM 3. Otherwise, we have
\[
\Phi(g(k)) \leq \Phi(g(l)) \leq n \ \mbox{and} \ g(k)!_n \subeq g(l)!_n,
\]
thus $g(k)\! \subeq \! g(l)$ by CLAIM 4. \qed

\section{Independence Results}\label{sec:ind}

Arguably, among the greatest milestones in proof theory are Gödel's incompleteness theorems \cite{Godel31} regarding the unprovability, or equivalently the independence, of some statements over a mathematical theory; in the case of the second incompleteness theorem the consistency of Peano arithmetic over Peano arithmetic itself. The first concrete example of independence (i.e., a true, yet unprovable, statement with a numerical, rather than logical, content) was found in 1977 by Paris and Harrington \cite{PH77}; since then, many other cases have been discovered \cite{KP82,Simpson85nonprovability,Smith85}.

Following the approach of \cite{FW23}, we consider tree-like structures proving that the well quasi-orderedness of some, relatively simple, embeddability relations on these structures is not provable over relevant theories, such as Peano Arithmetic, \textsf{PA}, or Arithmetic Transfinite Recursion \textsf{ATR$_0$}; see \cite{Simpson09} for a detailed presentation of \textsf{ATR$_0$}.

For what concerns notation let us define $2_0(l) \putAs l$ and $2_{k+1}(l) \putAs 2^{2_k(l)}$, so that $2_k(l)$ is an iterated exponential function of height $k$ with the number $l$ at the highest position.

\subsection{Independence Results for Ackermannian Terms}

In this section, we treat independence results, over first- and second-order arithmetical theories, for tree-like structures built upon Ackermannian terms. First, let us consider a suitable version of the Ackerman function.
\begin{defi}\label{def:Ack} Let $a,b,k$ be integers with $k \! \geq \! 2$, we define the following Ackermannian function $A_a(k,b) \! \in \! \mN$. Fix $k$ and let us denote $A_a(k,b)$ by $A_ab$ with $A_a(-1) \putAs 1$ as auxiliary value, then $A_ab$ is recursively defined as:\footnote{The iterations of a function $F$ are defined as: $F^0(i)\putAs i$ and $F^{n+1}(i) \putAs F(F^n(i))$. }
\[
\begin{array}{rcl} \vspace{0.1 cm}
A_0b & \putAs & k^b \\ \vspace{0.1 cm}
A_{a+1}b & \putAs & A_{a}^k A_{a+1}(b-1) \\ \vspace{0.1 cm}
A_{\omega}(k,b) & \putAs & A_bb \\
A_{\omega}(b) & \putAs & A_{\omega}(b,b).
\end{array}
\]
\end{defi}
Note that $A_0(k, 0) = 1$ regardless of $k$. Except for a few trivial cases, the Ackermann function increases strictly in each of its parameters. This can be confirmed by a straightforward inductive argument, which we leave to the reader.
\begin{lem} Let $a \leq a', b \leq b'$, and $2 \leq k \leq k'$ be natural numbers. Then,
\[
\max\{a, b\} < A_a(k, b) \leq A_{a'}(k',b').
\]
If, moreover, $a + b + k < a' + b' + k'$ and $a' + b' > 0$, then $A_a(k, b) < A_{a'}(k', b').$
\end{lem}
Exploiting these functions, for a fixed $k \geq 2$ the following set $T_k$ of Ackermannian terms is definable.
\begin{defi}\label{def:Tk} Let $T_k$ be the least set of terms such that:
\begin{enumerate}

\item $0 \in T_k$;

\item if $\mathfrak{a},\mathfrak{b} \in T_k \setminus \{0\}$, then $\mathfrak{a}+\mathfrak{b} \in T_k$;

\item if $\fa,\fb \in T_k$, then $A_\fa(k,\fb) \in T_k$.

\end{enumerate}
Moreover, let $T \putAs \bigcup_{k}^{\infty} T_k$ be the set of all Ackermannian terms.
\end{defi}
For $\fa \in T_k$, the \ita{value} of $\fa$, denoted by $|\fa| \! \in \! \mN$, is inductively defined as follows: $|0| \putAs 0$, $|\fa+\fb| \putAs |\fa| + |\fb|$, and $|A_\fa(k,\fb)| \putAs A_{|\fa|}(k,|\fb|)$. Furthermore, since $k$ is a parametric placeholder, given $k,h \geq 2$ a ``base change'' from elements of $T_k$ to elements of $T_h$ can be defined; more precisely, given $\fa \! \in \! T_k$, $\fa[k\putAs h]$ indicates the element of $T_h$ obtained from $\fa$ by replacing every occurrence of $k$ with $h$.

For each $k \geq 2$, the elements of $T_k$ can be represented as ordered binary trees with the leaves labelled by $0$ and the inner nodes labelled by $+$ or $A_k$. For example:
\[
A_{A_0(k,0)}(k,0) + A_0(k,A_0(k,0)) \  \  \ \leadsto \  \
\begin{tikzpicture}[baseline={([yshift=-.5ex]current bounding box.center)}, grow'=up, level distance=10mm,
level 1/.style={sibling distance=13mm},
level 2/.style={sibling distance=8mm}, 
level 3/.style={sibling distance=7mm},
 vertex/.style={anchor=base, circle,fill=black!25,minimum size=18pt,inner sep=2pt}]
    \usetikzlibrary {quotes}
       \node {$+$} [sibling distance=15mm]
        child {node {$A_k$}
                child {node {$A_k$}
                        child {node {$0$}}
                        child {node {$0$}}
                        }
                child {node {$0$}}
                }
        child {node {$A_k$}
                child {node {$0$}}
                child {node {$A_k$}
                        child {node {$0$}}
                        child {node {$0$}}
                      }
                };
  \end{tikzpicture}
\]
This connection allows to easily define the following relation on $T_k$.

\begin{defi}\label{def:<k} Let $\leq_k$ be the least binary relation on $T_k$ such that, for all $\fa,\fb,\fc \in T_k$, the following properties hold:
\begin{enumerate}

\item $ 0 \leq_k \fa$;

\item if $\fa \leq_k \fb$, then $\fa \leq_k \fb+\fc$ and $\fa \leq_k \fc+\fb$;

\item if $\fa \leq_k \fb$, then $\fa \leq_k A_\fb(k,\fc)$ and $\fa \leq_k A_\fc(k,\fb)$;

\item if $\fa \leq_k \fa'$ and $\fb \leq_k \fb'$, then $\fa+\fb \leq_k \fa' + \fb'$ and $A_\fa(k,\fb) \leq_k A_{\fa'}(k,\fb')$.

\end{enumerate}
\end{defi}
Let us observe that if $\fa \leq_k \fb$, then $|\fa| \leq |\fb|\! .$ Moreover, given the tree structure of its elements, $T_k$ is a proper subset of $\mT(\bar{3})$, with $\bar{3} \putAs\{0,+, A_k\}$, and it is easy to see that $\fa \leq_k \fb$ if and only if $\fa \preceq \fb$, where $\preceq$ is tree embeddability (see \Def{def:treemb}).

The next definition introduces the first principle under scrutiny.

\begin{defi}\label{def:W(f)} For a non-decreasing number-theoretic function $f$ of two arguments, let $W(f)$ be the following assertion:
\begin{multline*}
(\forall K)(\exists M)(\forall \fa_0, \dots, \fa_M \in T) \big[(\forall i \leq M)[\fa_i \! \in \! T_{i+2} \et |\fa_i| \leq f(K,i)] \\  \frec (\exists i,j \leq M) \left( i < j \, \et \, \fa_i[i+2 \putAs j+2] \leq_{j+2} \fa_j\right) \big].
\end{multline*}
\end{defi}
For what concerns the provability of $W(f)$ with respect to some specific functions $f$, we have the following two results.

\begin{theor}\label{theor:Sigmaf} Let $d \geq 1$ be a natural number and fix $f_d(K,i) \putAs 2_{d-1}((i+2)^K)$, then:

\begin{enumerate}

\item  $W(f_d)$ is true;

\item  $I\Sigma_d \, \nvdash \ W(f_d)$.

\end{enumerate}

\end{theor}


\begin{theor}\label{theor:ATRf} Let $f(K,i)\putAs A^K_{\omega}(i+2,0)$,\footnote{After the first binary application of $A_{\omega}(\sqdot,\sqdot)$, we apply the unary version of $A_{\omega}(\sqdot)$, as defined in the last line of \Def{def:Ack}.} then:

\begin{enumerate}

\item $W(f)$ is true;

\item  \textsf{ATR$_0$} $\, \nvdash \ W(f)$.

\end{enumerate}

\end{theor}
We consider, and prove, the following slight generalization. Denoted with $\mN_{\geq 2}$ the set of integers greater than $1$, let $t \!: \mN_{\geq 2} \times \mN \frec T$ be a selection function sending $(i,n) \mapsto t_i(n) \in T_i$ with the only constraint that $|t_i(n)|=n$. For any such selection function and any non decreasing number-theoretic function $f$ of two arguments, let $W(t, f)$ be the following
assertion:
\begin{multline*}
(\forall K)(\exists M)(\forall a_0, \dots, a_M \in \mN) \big[(\forall i \leq M)[ a_i \leq f(K,i)] \\  \frec (\exists i,j \leq M) \left( i < j \, \et \, t_{i+2}(a_i)[i+2 \putAs j+2] \leq_{j+2} t_{j+2}(a_j)\right) \big].
\end{multline*}
In particular, we consider the specific selection function that sends $(i,n)$ to the tree representation of the \ita{$i$-base Ackermannian normal form} of $n$. More precisely, as detailed in \cite{AFWW:unprovable}(see also \cite{AWW:parametrized,FW:walkG,FW:walkGA}), for a fixed $k \geq 2$, every natural number $n$ can be uniquely written in the form $A_a(k,b)+c$, for suitable natural numbers $a,b,c \geq 0$, if we impose some conditions (e.g., $a$ is the smallest number such that $A_{a+1}(k,b) > n$ and $b$ is the smallest one such that $A_a(k,b+1) > n$). Moreover, $a,b$ and $c$ can themselves be put into their $k$-base Ackermannian normal form and so on; we write $n =_k A_a(k,b)+c$ when $A_a(k,b)+c$ is the $k$-base Ackermannian normal form of $n$ and $a,b,c$ are themselves in their normal form. Since the $k$-base Ackermannian normal form involves only addition and the function $A_{\sqdot}(k,\sqdot)$, it can be written as an element of $T_k$; for example $0 \! =_2 \! 0, 1 \! =_2 \! A_00, 4 \! =_2 \! A_{A_00}0$ and $2^{16}+1 \! =_2 \! A_{A_00}A_00 + A_00$.

The aforementioned Ackermannian normal form and its associate selection function $t$ enjoy the following two combinatorial properties: the former is the so-called \ita{base-change maximal} property ensuring that the Ackermannian normal form maximizes the values with respect to change of base; the latter is a form of commutative between the Ackermannian normal form and the base-change operation.
\begin{lem}\label{lem:propNormAck} Let $t$ be the Ackermannian normal form selection function, then
\begin{itemize}

 \item $\forall k \! \geq \! 2 \ \forall \tau \! \in \! T_k\ \forall l \! > \! k \ \ |\tau[k:=l]| \, \leq \, \left|t_k(|\tau|)[k:=l]\right|$;
 
 \item $\forall k \! \geq \! 2 \ \forall \tau \! \in \! T_k\ \forall l \! > \! k \ \ t_k(|\tau|)[k:=l] \, = \, t_l \big( \left|t_k(|\tau|)[k:=l]\right| \big)$.

\end{itemize}
\end{lem}
\proof See \cite{AFWW:unprovable,FW:walkGA}. \qed

For such a selection function, we have the following independence results.
\begin{theor}\label{theor:Sigmaf} Let $d \! \geq \! 1$ be a natural number; if $f_d(K,i) \putAs 2_{d-1}((i+2)^K)$,~then:

\begin{enumerate}

\item  $W(t,f_d)$ is true;

\item  $I\Sigma_d \, \nvdash \ W(t,f_d)$.

\end{enumerate}

\end{theor}

\proof The proof is analogous to the one of the next theorem. \qed
\begin{theor}\label{theor:ATRf} Fix $f(K,i)\putAs A^K_{\omega}(i+2,0)$, then

\begin{enumerate}

\item $W(t,f)$ is true;

\item  \textsf{ATR$_0$} $\, \nvdash \ W(t,f)$.

\end{enumerate}

\end{theor}

\proof \ita{1.}(Correctness) For a given $K$, choose $M$ according to the following miniaturization of Kruskal's theorem
\begin{multline*}
(\forall K)(\exists M)(\forall t_0, \dots, t_M \in \mT(\bar{3}))\\
\left[(\forall i \leq M)[N(t_i) \leq 4 \cdot A^K_{\omega}(i+2,0) +1] \frec (\exists i < j \leq M)[ t_i \preceq t_j]\right]
\end{multline*}
where, similarly as before, $\mT(\bar{3})$ is the set of finite trees with labels in $\bar{3}=\{1,+,A_{x}\}$ and, for $t \! \in \! \mT(\bar{3})$, $N(t)$ is the number of nodes of $t$; see \cite{FW23,Gallier91} for a proof of this statement. Let $a_0, \dots, a_M$ be a given sequence such that $a_i \! \leq \! A^K_{\omega}(i+2,0)$, then $t_{i+2}(a_i)[i+2:=x] \! \in \! T_x \subseteq \mT(\bar{3})$ and $N(t_{i+2}(a_i)[i+2:=x]) \leq 4 \cdot A^K_{\omega}(i+2,0) +1$; the last inequality holds since an Ackermannian term $\fa$ with value $|\fa| \leq \! A^K_{\omega}(i+2,0)$, such as $t_{i+2}(a_i)$, has less then $4 \cdot A^K_{\omega}(i+2,0) +1$ nodes when seen as a tree. This non optimal upper bound derives from the fact that writing a number $n$ using just $1$ and $+$ required $2n-1$ symbols (i.e., nodes in a binary tree); in our case, each $1$ is written as $A_0(k,0)$ (two leaves labelled by $0$ and an inner node labelled by $A_k$) and thus requires in total $4n-1$ nodes. By the previous miniaturization of Kruskal's theorem, there exist $i<j \leq M$ such that $t_{i+2}(a_i)[i+2:=x] \preceq t_{j+2}(a_j)[j+2:=x]$, thus $t_{i+2}(a_i)[i+2 \putAs j+2] \leq_{j+2} t_{j+2}(a_j)$. Let us observe that an easy adaptation of this strategy can readily be used to prove point \ita{1.} of \Theor{theor:ATRf}.

\smallskip

\ita{2.}(Unprovability). For $2 \! \leq \! h \! \leq \! l$, let $a[h := l]$ be the result of replacing in the  $h$-base Ackermannian normal form of $a$ all occurrences of $h$ by $l$; namely, if $a =_h A_r(h,p) + q$, then $a[h := l] \putAs A_{r[h:=l]}(l, p[h:=l]) + q[h := l]$. The unprovability result relies on the following claim.

\noindent \textbf{Claim.} Let $a_0 :=A^K_{\omega}(2,0)$ and $a_{i+1}:=a_i[i+2 :=i+3]-1$, then:
\begin{enumerate}

 \item $\forall K\, \exists M'\ a_{M'}=0$ holds;
 
 \item \textsf{ATR$_0$} $\nvdash \, \forall K\, \exists M'\ a_{M'}=0$.
\end{enumerate}

\noindent \ita{Proof of the claim.} See \cite{AFWW:unprovable}.

Let $f(K,i)\! :=A^K_{\omega}(i+2,0)$ and, for a given $K$, choose $M$ according to $W(t,f)$. Moreover, as in the claim, let $a_0=A^K_{\omega}(2,0)$ and $a_{i+1}=a_i[i+2:=i+3]-1$; then, $a_i \! \leq \! A^K_{\omega}(i+2,0)$. Assume, by $W(t,f)$, that there exist $i \! < \! j \! \leq \! M$ such that $t_{i+2}(a_i)[i+2 \putAs j+2] \leq_{j+2} t_{j+2}(a_j)$; then, $|t_{i+2}(a_i)[i+2 \putAs j+2]| \leq \left|t_{j+2}(a_j)\right|=a_j$. Now we claim that, if $a_m \neq 0$ for all $0 \! \leq \! m \! \leq \! M$, then $|t_{i+2}(a_i)[i + 2 := j + 2]| \! > \! |t_{j+2}(a_j)|$. 
This contradiction shows that $W(t, f)$ implies $\forall K\, \exists M' \, a_{M'}=0$, thus \textsf{ATR$_0$} does not prove $W(t,f)$ for $f(K,i)\! :=A^K_{\omega}(i+2,0)$. 
The claim is proved by induction on $k \! > \! i$. Let $k=i+1$, then $|t_{i+2}(a_i)[i+2:=i+3]|=|t_{i+2}(a_i)|[i+2:=i+3]=a_i[i+2:=i+3]>a_i[i+2:=i+3]-1=a_{i+1}=|t_{i+3}(a_{i+1})|$; the first equality derives from the choice of $t$, i.e., $t_{i+2}(a_i)$ is the (tree representation of the) $i+2$-base Ackermannian normal form of $a_i$ (see \Lem{lem:propNormAck}). The induction step, $k \frec k+1$, is proved analogously, exploiting the fact that $a_k [k+2 := k+3] > a_k [k+2 := k+3]-1 = a_{k+1} = |t_{k+3}(a_{k+1})|$. \qed

\smallskip

Lastly, we prove the equivalence between $W(f)$ and $\forall t\, W(t,f)$; since, by the previous theorem \textsf{ATR$_0$} $\nvdash \, \forall t\, W(t,f)$, this implies that \textsf{ATR$_0$} $\nvdash W(f)$.
\begin{prop}[\textsf{ATR$_0$}] Let $f$ be a non-decreasing number-theoretic function of two arguments, then $W(f)$ and $\forall t\, W(t, f)$ are equivalent.
\end{prop}
\proof $``W(f) \imp \forall t\, W(t,f)$'' Let $t$ be any selection function, we need to prove that
\begin{multline*}
(\forall K)(\exists M)(\forall a_0, \dots, a_M \in \mN) \big[(\forall i \leq M)[ a_i \leq f(K,i)] \\  \frec (\exists i,j \leq M) \left( i < j \, \et \, t_{i+2}(a_i)[i+2 \putAs j+2] \leq_{j+2} t_{j+2}(a_j)\right) \big].
\end{multline*}
Let $K \! \in \! \mN$, $M$ be given as in $W(f)$ and consider a finite sequence $a_0, \dots, a_M\! \in \mN$ such that, for all $i \! \leq \! M$, $a_i \! \leq \! f(K,i)$; then, the finite sequence of Ackermannian terms $t_2(a_0), \dots, t_{M+2}(a_M)$ satisfies the hypotheses of $W(f)$, i.e., $t_{i+2}(a_i) \! \in \! T_{i+2} \et |t_{i+2}(a_i)|=a_i \leq f(K,i)$ for all $i \leq M$. Let $i<j \leq M$ as given by $W(f)$, we obtain $t_{i+2}(a_i)[i+2:=j+2] \leq_{j+2} t_{j+2}(a_j)$; thus $W(t,f)$ holds.

\smallskip

\noindent $``\forall t\, W(t,f) \imp W(f)$'' We reason by contradiction assuming $\neg W(f)$, namely
\begin{multline*}
(\exists \bar{K})(\forall M)(\exists \tau_0, \dots, \tau_M \in T) \big[(\forall i \leq M)[\tau_i \! \in \! T_{i+2} \et |\tau_i| \leq f(\bar{K},i)] \\  \et (\forall i <j \leq M) \left(\tau_i[i+2 \putAs j+2] \not\leq_{j+2} \tau_j\right) \big].
\end{multline*}
Thus, for every $M$, we can find a ``bad sequence'' $\tau_0, \dots, \tau_M$ such that, for all $i \leq M$, $\tau_i \in T_{i+2}$ and $|\tau_i| \leq f(\bar{K},i)$, but $\tau_i[i+2 \putAs j+2] \not\leq_{j+2} \tau_j$ for all $i<j \leq M$. Since such sequences are closed under initial segments and, for every $i$, there are finitely many terms $\tau$ in $T_{i+2}$ such that $|\tau| \leq f(\bar{K},i)$, this type of bad sequences form a finitely branching infinite tree. By König's lemma, we can obtain an infinite sequence $\tau_0, \tau_1, \dots$ such that $\tau_i \in T_{i+2}$ and $|\tau_i| \leq f(\bar{K},i)$ for all $i$, and $\tau_i[i+2 \putAs j+2] \not\leq_{j+2} \tau_j$ for all $i<j$. Given such sequence, we consider the following selection function $s$, where $t^{NF}$ is the Ackermannian normal form selection function.
\[
s(i+2,n) \putAs \left\{
 \begin{array}{ll}\vspace{0.05 cm}
   \tau_i & \mbox{if}\, n=|\tau_i|, \\
   t_{i+2}^{NF}(n) & \mbox{otherwise}.
 \end{array}
 \right.
\]
From $W(s,f)$, with respect to $\bar{K}$, we obtain a value $M_{\bar{K}}$ such that
\begin{multline*}
(\forall a_0, \dots, a_{M_{\bar{K}}} \in \mN) \big[(\forall i  \leq  M_{\bar{K}})[ a_i \leq f(\bar{K},i)] \\ \frec (\exists i <  j  \leq  M) \left( s_{i+2}(a_i)[i+2 \putAs j+2] \leq_{j+2} s_{j+2}(a_j)\right) \big].
\end{multline*}
Applying this property to the sequence $|\tau_0|, |\tau_1|, \dots, |\tau_{M_{\bar{K}}}|$, there are indexes $i < j \leq M_{\bar{K}}$ such that
\[
\tau_i[i+2:=j+2]=s_{i+2}(|\tau_i|)[i+2:=j+2] \leq_{j+2} s_{j+2}(|\tau_j|) = \tau_j
\]
which contradicts our hypothesis. \qed

\subsection{Independence Results for Exponential Expressions}

In this paragraph, we treat an independence result over Peano arithmetic regarding exponential expressions.
\begin{defi}\label{def:EXP} Let $EXP$ the least set of terms in one variable $x$ such that:
\begin{enumerate}

\item $0 \in EXP$;

\item if $\fa,\fb \in EXP$, then $x^\fa + \fb \in EXP$.

\end{enumerate}
\end{defi}
For $\fa \in EXP$ and a positive integer $k$, we denote by $\fa(k)$ the integer value of the term $\fa$ when the variable $x$ is set to $k$; namely, $\fa(k)$ is inductively defined as: $\fa(k)\putAs 0$ if $\fa \equiv 0$, and $\fa(k)\putAs k^{\fc(k)} + \fb(k)$ if $\fa \equiv x^\fc +\fb$. Moreover, $\fa(\omega)$ indicates the ordinal, smaller than $\ep_0$, obtained by replacing $x$ with $\omega$ and interpreting $+$ as the natural sum between ordinals (see \cite{Bachmann55} for the definition of natural sum). Lastly, similar to how we defined the exponential iteration with base $2$, let us define $\omega_0 \putAs 1$ and $\omega_{k+1} \putAs \omega^{\omega_k}$, so that $\omega_k$ is an $\omega$-tower of height $k$.

Our goal is to study a Slow Well-Ordering principle.

\begin{defi}\label{def:SWO} Let $SWO$ be the following assertion:
\begin{multline*}
(\forall K)(\exists M)(\forall \fa_0, \dots, \fa_M \in EXP)\\
\left[(\forall i \leq M)[\fa_i(2) \leq 2_K(i)] \frec (\exists j < M)[\fa_j(2_K(M)) \leq \fa_{j+1}(2_K(M))]\right].
\end{multline*}
\end{defi}
For what concerns the provability of $SWO$, we obtain the next independence result.
\begin{theor}\label{theor:SWO} \phantom{placeholder}
\begin{enumerate}

\item $SWO$ is true;

\item $\textsf{PA} \nvdash \, SWO$.

\end{enumerate}

\end{theor}

\proof \ita{1.}(Correctness) The truth of $SWO$ follows from the corresponding result for $\ep_0$ with respect to the maximal coefficient. Namely, we consider the following statement:
\begin{multline}\label{eq:SWO}
(\forall K)(\exists M)(\forall \al_0, \dots, \al_M \leq \ep_0)\\
\big[(\forall i  \leq  M)[|\al_i|\leq 2_K(i)] \frec (\exists j  <  M)[\al_j \leq \al_{j+1}]\big]
\end{multline}
where, for an ordinal $\al < \ep$ in Cantor normal form (i.e., $\al = \omega^{\al_1} \cdot a_1 + \dots + \omega^{\al_r} \cdot a_r$ with $\ep_0 \! > \! \al_1 \! > \! \dots \! > \! \al_r$ and $a_1, \dots, a_r \! \in \! \mN_{+}$), the \ita{maximal coefficient} $|\al|$ is recursively defined as $|\al| \putAs \max\{|\al_1|, \dots, |\al_r|, a_1, \dots, a_r\}$; see \cite{Smith85} for the proof of \eqref{eq:SWO}. 
Let $K$ be given and choose $M$ with respect to \eqref{eq:SWO}. If $\fa_0, \dots , \fa_M$ is a sequence in $EXP$ such that $\fa_i(2) \leq 2_K(i)$ for all $i \leq M$, then $\fa_i(\omega)$ is a sequence of ordinals less than $\ep_0$ such that the maximal coefficient of $\fa_i(\omega)$ is bounded by $\fa_i(2)$, hence by $2_K(i)$. The assertion for the maximal coefficient yields the existence of a natural number $j < M$ such that $\fa_j(\omega) \leq \fa_{j+1}(\omega)$. We have $\fa_i(2) \leq 2_K(i)$ for all $i \leq M$, therefore the maximal coefficient of $\fa_j(\omega)$ is strictly smaller than $2_K(M)$. 
For the last step of correctness, we rely on some properties of the \ita{slow-growing hierarchy} $G$ that, up to $\ep_0$, is inductively defined as: $G_{0}(x) \putAs 0, G_{\al+1}(x) \putAs G_{\al}(x) +1$, and $G_{\la}(x) \putAs G_{\la[x]}(x)$ if $\la$ is a limit ordinal.\footnote{Here and below, the standard assignment of fundamental sequences to ordinals smaller that $\ep_0$ is used; see \cite[Definition 3]{MW15} or \cite{BCW94}.} 
In our case, we find $G_{\fa_j(\omega)}(2_K(M)) \leq G_{\fa_{j+1}(\omega)}(2_K(M))$ and, since $G_{\fa_j(\omega)}(2_K(M)) = \fa_j(2_K(M))$ and $G_{\fa_{j+1}(\omega)}(2_K(M)) = \fa_{j+1}(2_K(M))$, the conclusion follows. The two required properties of the slow-growing hierarchy $G$, namely $\al \! \leq \! \be, |\al| \! \leq \! x \Rightarrow G_{\al}(x) \leq G_{\be}(x)$ and $G_{\al}(x) = \al[\omega \putAs x]$, are proved in~\cite[Lemma 6]{MW15}.

\ita{2.}(Unprovability) Assume $K \geq 2$. Let us define $\al_0 \putAs \omega_{K-2}$ and $\al_{i+1} \putAs \al_i[i]$, and let $N$ be the least $i$ such that $\al_i = 0$; then \textsf{PA} does not prove $\forall K\, \exists N\, \al_N \! = \! 0$ \cite{BCW94}. The maximal coefficient of $\al_i$ is smaller than $i$ for all $0 < i \leq M$ and, since $\al_i < \al_0$, we find $G_{\al_i}(i) < G_{\omega_K}(i) \leq 2_K(i)$.

Put $\fa_i \putAs \al_i[\omega := x]$, i.e., any occurrence of $\omega$ in $\al_i$ is replaced by $x$; then $\fa_i \in EXP$. Moreover, for $i \geq 2$ we obtain $\fa_i(2) \leq \al_i[\omega := 2] = G_{\al_i}(2) \leq G_{\al_i}(i) \leq G_{\omega_K}(i) \leq 2_K(i)$; while, for $i < 2$, $\fa_i(2) \leq 2_K(i)$ can be found by direct evaluation. Let us apply $SWO$ to $K$ in order to obtain an appropriate $M$, and assume, to reach a contradiction, that $\al_M \neq 0$, i.e., $\al_i \! > \! \al_{i+1}$ for all $i < M$. If we consider the sequence $\fa_i$, on the one hand we obtain, by $SWO$, that there exists an $j < M$ such that $\fa_j(2_K(M)) \leq \fa_{j+1}(2_K(M))$. On the other hand, we have that $\al_i\! >\! \al_{i+1}$, and the maximal coefficient of $\al_i$ is strictly bounded by $2_K(M)$ which leads to $G_{\al_i}(2_K(M)) > G_{\al_{i+1}}(2_K(M))$. Hence $\fa_j(2_K(M)) > \fa_{j+1}(2_K(M))$, contradiction. This shows that $SWO$ implies $\forall K\, \exists N\, \al_N =0$ and thus, \textsf{PA} does not prove $SWO$. \qed

\section{Related and Future Work}\label{sec:work}

Well quasi-orders have been extensively studied in the recent decades, with fruitful applications in algebra \cite{AP03}, graph theory \cite{Pouzet85}, logic (remarkably proof theory \cite{Simpson85nonprovability}) and even theoretic computer science (e.g. for term rewriting \cite{Dershowitz82,Dershowitz87,DO88} and program termination \cite{BG08}). For a general introduction and an historical survey, we refer respectively to \cite{Milner85} and \cite{Kruskal72}; whereas for a recent collection that illuminates the multifaceted nature of {\wqo}, we refer to \cite{SSW20}.

For what concerns Higman's lemma, recent investigations, particularly in reverse mathematics, are related to its extension called ``generalized Higman's lemma''. As previously said, well quasi-orders are not preserved under infinitary operations, such as infinite sequences; to overcome this obstacle better quasi-orders, bqo, have been introduced \cite{Nash-W68}. Involving bqo, the generalized Higman's lemma reads as follows: ``if $Q$ is a bqo, then $Q^*$ is a bqo''. It can be proved (see for example \cite[Theorem 6.21]{Marcone20}) that the generalized Higman's lemma is actually equivalent to the celebrated Nash-Williams' theorem \cite{Nash-W68}: denoted with $\tilde{Q}$ the set of countable sequences over $Q$, i.e., transfinite sequences of countable length, Nash-Williams' theorem states that ``if $Q$ is a bqo, then $\tilde{Q}$ is a bqo''.

Regarding Kruskal's theorem, the latest developments, attributed to Anton Freund, Patrick Uftring, and collaborators \cite{Freund20gap,Freund:uniform,FRW22,FU23,Uftring2025}, primarily focus on the so-called ``uniform Kruskal theorem". This theorem extends the classical result for trees to encompass general recursive data types and incorporates Girard's dilator theory \cite{Girard81,Girard85}, giving it a pronounced categorical character. It is worth mentioning also the rich literature related to constructive proofs of Higman's lemma \cite{CF93,MR90,RS93,Seisen01Higman} and Kruskal's theorem \cite{Seisen01Kruskal,Veldman04}. The most recent result in this stream is a novel version of Higman's lemma relative to bars \cite{BBS24}.

For what concerns independence results, the ancestors are obviously Gödel's incompleteness theorems \cite{Godel31}. Despite the fundamental theoretical relevance of Gödel's achievement, the first example of a concrete mathematical statement independent from PA was only discovered in 1977 by Paris and Harrington \cite{PH77}. Following their finding, many other independence results have been discovered \cite{KP82,Simpson85nonprovability,Smith85}. In more recent years, the study of independence results exposes the so-called ``phase transition'' phenomenon in logic \cite{GW10,Weiermann05} where the shift from provability to unprovability is characterized by a sharp threshold.

\smallskip

With respect to future developments regarding Kruskal's theorem, the proof-theoretic work of Rathjen and Weiermann \cite{RW93} continues to serve as a guiding reference. Beyond calculating its proof ordinal, namely \RCA\, $\vdash \mbox{KT}(\omega) \sse \mbox{WO}\left( \vartheta (\Omega^{\omega}) \right)$, they characterized the proof-theoretic strength of Kruskal's theorem for unlabelled trees in terms of reflection principles.

\begin{theor}[Rathjen and Weiermann] 
\[
\mbox{RCA}_0 \vdash \, \mbox{KT}(\omega) \sse \Pi^1_1\mbox{-}RF\!N(\Pi^1_2\mbox{-}BI_0).
\]

\end{theor}
Where $\Pi^1_1\mbox{-}RFN(\Pi^1_2\mbox{-}BI_0)$ stands for the \ita{uniform reflection principle} for $\Pi^1_1$ formulas of the theory $\Pi^1_2\mbox{-}BI_0$, which amounts to \RCA\, extended with bar induction for $\Pi^1_2$ formulas; see \cite[Section 11]{RW93} for more details. 

As a result, the next step is to achieve a comparable classification for \KTlw \ as well as $\forall n\,$\KTln. Currently, we have the following two conjectures:
\begin{conj}\phantom{placeholder}

\begin{itemize}

\item \RCA\, $\vdash \, \mbox{KT}_{\ell}(\omega) \sse \Pi^1_2\mbox{-}\omega RFN(\Pi^1_2\mbox{-}BI_0 \upharpoonright \Pi^1_3)$ \ \ [by F. Pakhomov]

\item \RCA\, $\vdash \, \forall n\, \mbox{KT}_{\ell}(n) \sse \Pi^1_2\mbox{-}RFN(\Pi^1_2\mbox{-}BI_0)$ \ \ \ \ \ \ [by A. Freund].

\end{itemize}

\end{conj}
Further research in this area will be needed to expand the current results to encompass other embeddability relationships between trees, such as Friedman’s gap condition \cite{Freund20gap,Kriz89}.

For what concerns independence results, if terms for the Ackermann function are added to $EXP$, we expect to achieve an independence result for \textsf{ATR$_0$}; moreover, integrating terms for the fast-growing hierarchy $(F_{\al})_{\al <\ep_0}$ should lead to independence results for $ID_1$.

\section*{Conclusion}\label{sec:conclusion}

This paper has treated the proof-theoretic relations between Higman's and Kruskal's theorems in the framework of reverse mathematics considering both their finite and infinite versions. 
After a brief summary of {\wqo} theory in reverse mathematics, including also novel results regarding the closure of {\wqo} under order-preserving and order-reflecting functions, the main elements of Higman's and Kruskal's theorems, namely ordered algebras and labelled trees, have been introduced. The equivalence between the finite versions of the two theorems, denoted respectively as $\forall n$ HT$(n)$ and $\forall n$ \KTln, has been thoroughly presented, followed by the salient points for the infinitary case, i.e., HT$(\omega)$ and \KTlw. The previous equivalences, together with some other recent findings related to Kruskal's theorem, have been summarized in an exhaustive schema depicting the proof-theoretic relations over \RCA\, of the results under consideration. 
Furthermore, also the connections between Higman's lemma and Dickson's lemma have been scrutinized, collecting the main references concerning their relations and proof-theoretic ordinals. 

In addition, and enriching previous findings, some independence results over first- and second-order theories were treated. In particular, tree-like structures, involving either Ackermannian terms or exponential expressions, have been studied unveiling well-foundedness properties that are independent from Peano arithmetic and relevant fragments of second-order arithmetic.

Finally, two conjectures regarding the proof-theoretic strength of the aforementioned versions of Kruskal's theorems with labels have been stated.

\subsection*{Acknowledgments}

The authors thank Peter Schuster and Anton Freund for their valuable suggestion during the preparation of this paper, which benefited from the fruitful discussions during the 2025 meeting on Reverse Mathematics at the Erwin Schrödinger Institute in Vienna. The authors thank also the anonymous referee for their thorough reading and beneficial proposals.

The first author is member of the “Gruppo Nazionale per le Strutture Algebriche, Geometriche e le loro Applicazioni” (GNSAGA) of the ``Istituto Nazionale di Alta Matematica'' (INdAM).

No use of AI was made during the writing of this paper.

\bibliographystyle{plain} 
\bibliography{HigmanBiBold}  

\end{document}